\documentclass[12pt]{article}
\usepackage{amssymb,amsmath,amscd, amsthm,stmaryrd,verbatim, mathrsfs}
\numberwithin{equation}{section}

\begin{document}
\baselineskip=14pt

\newcommand{\la}{\langle}
\newcommand{\ra}{\rangle}
\newcommand{\psp}{\vspace{0.4cm}}
\newcommand{\pse}{\vspace{0.2cm}}
\newcommand{\ptl}{\partial}
\newcommand{\dlt}{\delta}
\newcommand{\sgm}{\sigma}
\newcommand{\al}{\alpha}
\newcommand{\be}{\beta}
\newcommand{\G}{\Gamma}
\newcommand{\gm}{\gamma}
\newcommand{\vs}{\varsigma}
\newcommand{\Lmd}{\Lambda}
\newcommand{\lmd}{\lambda}
\newcommand{\td}{\tilde}
\newcommand{\vf}{\varphi}
\newcommand{\yt}{Y^{\nu}}
\newcommand{\wt}{\mbox{wt}\:}
\newcommand{\rd}{\mbox{Res}}
\newcommand{\ad}{\mbox{ad}}
\newcommand{\stl}{\stackrel}
\newcommand{\ol}{\overline}
\newcommand{\ul}{\underline}
\newcommand{\es}{\epsilon}
\newcommand{\dmd}{\diamond}
\newcommand{\clt}{\clubsuit}
\newcommand{\vt}{\vartheta}
\newcommand{\ves}{\varepsilon}
\newcommand{\dg}{\dagger}
\newcommand{\tr}{\mbox{Tr}}
\newcommand{\ga}{{\cal G}({\cal A})}
\newcommand{\hga}{\hat{\cal G}({\cal A})}
\newcommand{\Edo}{\mbox{End}\:}
\newcommand{\for}{\mbox{for}}
\newcommand{\kn}{\mbox{ker}}
\newcommand{\Dlt}{\Delta}
\newcommand{\rad}{\mbox{Rad}}
\newcommand{\rta}{\rightarrow}
\newcommand{\mbb}{\mathbb}
\newcommand{\lra}{\Longrightarrow}
\newcommand{\X}{{\cal X}}
\newcommand{\Y}{{\cal Y}}
\newcommand{\Z}{{\cal Z}}
\newcommand{\U}{{\cal U}}
\newcommand{\V}{{\cal V}}
\newcommand{\W}{{\cal W}}
\newcommand{\sta}{\theta}
\setlength{\unitlength}{3pt}
\newcommand{\msr}{\mathscr}
\newcommand{\wht}{\widehat}

\begin{center}{\large \bf Full Projective Oscillator Representations of Special \\ \pse Linear Lie Algebras and Combinatorial Identities} 
\end{center}
\vspace{0.2cm}

\begin{center}{\large Zhenyu Zhou\footnote{Corresponding author.} and Xiaoping Xu\footnote{Research supported
 by National Key R\&D Program of China 2020YFA0712600.
} }\end{center}

\begin{center}{
HLM, Institute of Mathematics, Academy of Mathematics \& System
Sciences\\ Chinese Academy of Sciences, Beijing 100190, P.R. China
\\ \& School of Mathematics, University of Chinese Academy of Sciences,\\ Beijing 100049, P.R. China}\end{center}

\begin {abstract}
\quad

Using the projective oscillator representation of $sl(n+1)$ and Shen's mixed
product for Witt algebras, Zhao and the second author (2011) constructed a new functor from
$sl(n)$-{\bf Mod} to $sl(n+1)$-{\bf Mod}. In this paper, we start from $n=2$ and use the functor successively to obtain a full projective oscillator realization of any  finite-dimensional irreducible representation of $sl(n+1)$ . The representation formulas of all the root vectors of $sl(n+1)$ are given in terms of first-order differential operators in $n(n+1)/2$ variables. One can use the result to study tensor decompositions of finite-dimensional irreducible modules by solving certain first-order linear partial differential equations, and thereby obtain the corresponding physically interested Clebsch-Gordan coefficients and  exact solutions of Knizhnik-Zamolodchikov equation in WZW model of conformal field theory.

 \vspace{0.3cm}

\noindent{\it Keywords}:\hspace{0.3cm} special linear Lie
algebra; projective oscillator
 representation; irreducible module;  singular vectors; combinatorial identities.

\end{abstract}

\section {Introduction}

 Sum-product type identities are important objects both in combinatorics and number theory. The Jacobi triple identity and
quintuple product identity are such well-known examples. Macdonald \cite{Mi} used affine analogues of the root systems of finite-dimensional simple Lie algebras to derive new type sum-product identities with the above two identities as special cases. Kac \cite{K} found the character formula for the integrable modules of affine
Kac-Moody algebras and showed that  the Macdonald identities are exactly the denominator identities. Kang and Kim \cite{Ks} found a number of interesting combinatorial identities from various expressions of the Witt partition functions in connection with the denominator identity of a certain graded Lie algebra.
The denominator identities of finite-dimensional simple Lie algebras are Vandermonde determinant type identities, which do not produce sum-product identities of numbers. In this paper, we show that our full projective oscillator representations of special linear Lie algebras naturally give rise to certain sum-product type identities of finite type.

Finite-dimensional irreducible representations of finite-dimensional simple Lie algebras over $\mbb C$ were abstractly determined
by Cartan and Weyl in early last century. However, explicit representation formulas of the root vectors in these simple algebras are in general very difficult to be given. In 1950, Gelfand and Tsetlin \cite{GT1, GT2}
used a sequence of corank-one subalgebras to obtain a basis whose elements were labeled by upside-down triangular data for finite-dimensional irreducible representations of general linear Lie algebras and orthogonal Lie algebras, respectively. Moreover, the corresponding matrix elements of the actions of simple root vectors (or Chevalley basis) were explicitly given.
Based on the theory of Mickelsson algebras and the representation theory of the Yangians, Molev \cite{M1} constructed a weight basis for finite-dimensional irreducible representations of symplectic Lie algebras and obtained explicit formulas for the matrix elements of generators. He \cite{M2} has also done similar work for $o(2n+1)$. There are also many interesting works  in this direction (e.g., cf. \cite{GK,G3,Lp,M1,FGR1,FGR2}). Knowing the representation formulas of simple roots is not enough to solve the general decomposition problem of
the tensors of irreducible representations because they are not commuting operators.

In this paper, we present a first-order differential operator realization of any finite-dimensional irreducible representation
of $sl(n+1)$ in $n(n+1)/2$ variables (cf. Theorem 2.3). Moreover, the explicit formulas of all the root vectors are given, which will be helpful in solving the general decomposition problem of the tensors of irreducible representations.
 In physics, the Clebsch-Gordan  coefficients are numbers that arise in angular momentum coupling in quantum mechanics. They appeared as the expansion
coefficients of total angular momentum eigenstates in an uncoupled tensor product basis (e.g., cf. \cite{E, dS}). We refer \cite{A, B,L,R} for more applications and later developments. Mathematically, the numbers are those of explicitly determining the irreducible components in the tensor of two finite-dimensional irreducible modules in terms of orthonormal bases. The first fundamental step is to find explicit formulas for the highest-weight vectors of those irreducible components. Even that is in general a very difficult problem.
 Our result in this paper simplifies the problem to solving certain first-order linear partial differential equations.
As a very special example, we find the explicit formulas of the highest-weight vectors of irreducible components in the tensor module of the irreducible module with highest weight $k\lmd_1$ and any finite-dimensional irreducible module of special linear Lie algebras. The well-known $sl(n+1)\downarrow sl(n)$ branching rule and the multiplicity-one theorem of $gl(n+1) \downarrow gl(n)$ are direct consequences. Our representation formulas can also be used to find exact solutions of Knizhnik-Zamolodchikov equation in WZW model of conformal field theory (cf. \cite{KZ, TK, X1}). An important feature of our representation is its connection with certain sum-product combinatorial identities of finite type. Below we give a more detailed introduction.

Let $n>1$ be an integer. Denote by $GL(n+1,\mbb{R})$ the group of $(n+1)\times (n+1)$
invertible matrices.  A {\it projective transformation}
 on $\mathbb{R}^n$ is given by
\begin{equation}u\mapsto \frac{Au+\vec b}{\vec c\:^t u+d}\qquad\mbox{for}\;\;u\in
\mathbb{R}^n,\end{equation} where all the vector in $\mathbb{R}^n$
are in column form and
\begin{equation}\left(\begin{array}{cc}A&\vec b\\ \vec c\:^t&
d\end{array}\right)\in GL(n+1,\mbb{R}).\end{equation} It is
well-known that a transformation of mapping straight lines to straight lines
must be a projective transformation. The group of projective
transformations is the fundamental group of $n$-dimensional
projective geometry.  Physically, the group  with $n=4$ consists of all the transformations of
keeping free particles including light signals moving with constant
velocities along straight lines (e.g., cf. \cite{GHW1, GHW2}). Based on the
embeddings of the Poincar\'{e} group and De Sitter group into the
projective group with $n=4$ , Guo, Huang and
Wu \cite{GHW1, GHW2} proposed three kinds of special relativity.

For simplicity, we assume that the base field is the field $\mbb C$ of complex numbers in the rest of this paper.
Let $E_{r,s}$ be the square matrix with 1 as its $(r,s)$-entry and 0
as the others. The special linear Lie algebra
\begin{equation}
sl(n+1)=\sum_{1\leq i<j\leq
n+1}(\mbb{C}E_{i,j}+\mbb{C}E_{j,i})+\sum_{r=1}^n\mbb{C}(E_{r,r}-E_{r+1,r+1}).\end{equation}
Let $\msr A=\mbb C[x_1,x_2,...,x_n]$ be the polynomial algebra in $n$ variables. Set
\begin{equation}D=\sum_{s=1}^nx_s\ptl_{x_s},\quad I_n=\sum_{i=1}^nE_{i,i}\in gl(n). \end{equation}
Let $M$ be an $sl(n)$-module. We fix $c\in \mbb C$ and make $M$ a $gl(n)$-module by letting $I_n|_M=c\:\mbox{Id}_M$. Denote
\begin{equation}\wht M={\msr A}\otimes_{\mbb C} M. \end{equation}
For any
two integers $p\leq q$, we denote $\ol{p,q}=\{p,p+1,\cdots,q\}$. Differentiating the transformations in (1.1), we get
 an inhomogeneous first-order differential operator representations of $sl(n+1)$. Using this representation and
  Shen's mixed product for Witt algebras in \cite{S}, Zhao and the second author
 \cite{ZX} obtained the following representation of $sl(n+1)$ on $\wht M$:
\begin{equation}E_{i,j}|_{\wht M}=x_i\ptl_{x_j}\otimes \mbox{Id}_M+\mbox{Id}_{\msr A}\otimes E_{i,j}|_M,\end{equation}
\begin{equation}(E_{i,i}-E_{j,j})|_{\wht M}=(x_i\ptl_{x_i}-x_j\ptl_{x_j})\otimes \mbox{Id}_M+\mbox{Id}_{\msr A}\otimes (E_{i,i}-E_{j,j})|_M,\end{equation}
\begin{equation}E_{i,n+1}|_{\wht M}=x_iD\otimes \mbox{Id}_M+x_i\otimes I_n|_M+\sum_{r=1}^nx_r\otimes E_{i,r}|_M,\;\;E_{n+1,i}|_{\wht M}=-\ptl_{x_i}\otimes \mbox{Id}_M,\end{equation}
\begin{equation}(E_{n,n}-E_{n+1,n+1})|_{\wht M}=(D+x_n\ptl_{x_n})\otimes \mbox{Id}_M+\mbox{Id}_{\msr A}\otimes (I_n+E_{n,n})|_M\end{equation}
 for $i,j\in\ol{1,n}$ with $i\neq j$. The first factors in the first terms of (1.6)-(1.9) are from differentiating (1.1). Moreover,\begin{equation}\ol{M}=U(sl(n+1))(1\otimes M)\end{equation} forms an $sl(n+1)$-submodule.   Denote by $\msr A_k$ the subspace of polynomials in $\msr A$ with degree $k$. Set
  \begin{equation}\wht M_k={\msr A}_k\otimes_{\mbb C} M,\quad \ol{M}_k=\ol{M}\bigcap\wht M_k. \end{equation}
 According to (1.6) and (1.7), $\wht M_k$ is the tensor module of the $sl(n)$-modules ${\msr A}_k$ and $M$, and $\ol{M}_k$ is an
$sl(n)$-submodule of $\wht M_k$.
 If $M$ is an irreducible $sl(n)$-module, then $\ol{M}$ is an irreducible $sl(n+1)$-module (cf. \cite{ZX, X2}).
 When $M$ is a finite-dimensional irreducible $sl(n)$-module, Zhao and the second author \cite{ZX} found a sufficient condition for $\wht M$ to be an irreducible $sl(n+1)$-module; equivalently $\wht{M}=\ol M$.

  When $M$ is a  highest-weight irreducible $sl(n)$-module,
 \begin{equation}\ol{M}=V_n(\lmd)\end{equation}is an irreducible highest-weight $sl(n+1)$-module for some highest weight $\lmd$.
 As we will see that $\lmd$ can be any dominant integral weight of $sl(n+1)$. Starting from $n=2$ and a finite-dimensional irreducible first-order differential operator representation of $sl(2)$, we apply (1.5)-(1.9) inductively in this paper to obtain a first-order differential operator realization of any finite-dimensional irreducible representation
of $sl(n+1)$ in $n(n+1)/2$ variables. Denote by  $\mbb N$ the set of nonnegative integers. Suppose that $\lmd=\sum_{i=1}^nk_i\lmd_i$ is a dominant integral weight of $sl(n+1)$.
We find that the degree of $\ol{M}$:
\begin{equation}\max\{r\in\mbb N\mid \ol{M}_r\neq\{0\}\}=\sum_{i=1}^nk_i=|\lmd|. \end{equation}
Moreover, $sl(n)$-modules
\begin{equation}\ol{M}_0=M=V_{n-1}(\sum_{i=1}^{n-1}k_{i+1}\lmd_i),\;\;
\ol{M}_{|\lmd|}=V_{n-1}(\sum_{i=1}^{n-1}k_i\lmd_i).\end{equation}
Indeed, $sl(n)$-modules $\ol{M}_r$ and $\ol{M}_{|\lmd|-r}$ have the same number of irreducible components.

According to Weyl's character formula,
\begin{equation}d_n(\lmd)=\dim V_n(\lmd)=\prod_{1\leq i<j\leq
n+1}\frac{\sum_{r=i}^{j-1}k_r+j-i}{j-i}\end{equation}
(e.g., cf. \cite{H, X2}). If $\lmd=k\lmd_i$, then the equation $\sum\limits_{r=0}^{|\lmd|}\dim \ol M_r=d_n(\lmd)$
is directly equivalent to the well-known classical combinatorial identity
    \begin{equation}
 \sum\limits_{r=0}^{k}\binom{k-r+i-1}{i-1}\binom{r+n-i}{n-i}=\binom{k+n}{k}\end{equation}
(e.g., cf. Page 10 in \cite{Ri}). As a byproduct, we prove
\begin{equation}\sum\limits_{s=0}^{n-1}(-1)^{s}\binom{\sum\limits_{j=s+1}^n(k_j+1)}{n}d_{n-1}(\sum\limits_{i=1}^sk_i\lambda_i+\lambda_s+\sum\limits_{i=s}^{n-1}k_{i+1}\lambda_i) =d_n(\sum_{\ell=1}^nk_\ell\lmd_\ell).\end{equation}
When $k_1=k_2=\cdots=k_n=k$, $\ol{M}$ is a Steinberg module and the above equation is equivalent to another well-known classical combinatorial identity
\begin{equation}
	\sum\limits_{s=0}^{n-1}(-1)^{s}\binom{(n-s)(k+1)}{n}\binom{n}{s}=(k+1)^n
\end{equation}	
(e.g., cf. Page 51 in \cite{Ri}).

\pse

In Section 2, we present the inductive construction of the differential-operator realization of any finite-dimensional
representation of $sl(n+1)$. In Section 3, we determine a basis for any finite-dimensional
representation of $sl(n+1)$ and its relation with Gelfand-Tzetlin basis. Moreover, we find all the  highest-weight vectors
for the irreducible $sl(n)$-components in $\wht M$ and $\ol M=V_n(\lmd)$. Section 4 is devoted to the detailed study on the special cases when the highest weights are
$k_1\lmd_1+k_2\lmd_2, k_1\lmd_1+k_2\lmd_n$ and $k\lmd_i$, respectively.

\section{General Construction}

In this section, we start from $n=2$ and repeatedly use (1.5)-(1.9) to construct the objective representation of $sl(n+1)$.
We also use the fact that the tensor algebra of two polynomial algebras is isomorphic to the polynomial algebra in all involved variables.

Given any positive integer $n>1$. Recall the special linear Lie algebra $sl(n+1)$ in (1.3).
Set
\begin{equation} h_i=E_{n+2-i,n+2-i}-E_{n+1-i,n+1-i},\qquad i\in\ol{1,n}. \end{equation}
The subspace
\begin{equation} H=\sum_{i=1}^n\mbb Ch_i\end{equation}
forms a Cartan subalgebra of $sl(n+1)$. We choose
\begin{equation} \{E_{i,j}\mid 1\leq j<i\leq n+1\}\;\;\mbox{as positive root vectors}.\end{equation}
In particular, we have
\begin{equation} \{E_{n+2-i,n+1-i}\mid i=1,2,...,n\}\;\;\mbox{as positive simple root vectors}.\end{equation}
Accordingly,
\begin{equation} \{E_{i,j}\mid 1\leq i<j\leq n+1\}\;\;\mbox{are negative root vectors}\end{equation}
and we have
\begin{equation} \{E_{n+1-i,n+2-i}\mid i=1,2,...,n\}\;\;\mbox{as negative simple root vectors}.\end{equation}
The representation formulas in \cite{GT1} were given only for the elements in (2.1), (2.4) and (2.6).
In particular,\begin{equation}
	sl(n+1)_-=\sum_{1\leq i<j\leq
		n+1}\mbb CE_{i,j}\;\;\mbox{and}\;\;sl(n+1)_+=\sum_{1\leq i<j\leq
		n+1}\mbb CE_{j,i}\end{equation} are the nilpotent subalgebra of
negative root vectors and the nilpotent subalgebra of positive root
vectors, respectively.
A  {\it singular vector} of an $sl(n+1)$-module $V$ is a weight vector
annihilated by the elements all positive root vectors. The fundamental
weights $\lmd_i\in H^\ast$ are
\begin{equation} \lmd_i(h_r)=\dlt_{i,r}.\end{equation}

Set
\begin{equation}\msr G_0=sl(n)+\mbb C(E_{n,n}-E_{n+1,n+1}),\;\;\msr G_+=\sum_{i=1}^n\mbb CE_{n+1,i},\;\msr G_-=
	\sum_{j=1}^n\mbb CE_{j,n+1}.\end{equation}
Then $\msr G_0$ is a Lie subalgebra of $sl(n+1)$ and $\msr G_\pm$ are abelian Lie subalgebras of $sl(n+1)$. Moreover,
\begin{equation}sl(n+1)=\msr G_-\oplus\msr G_0\oplus\msr G_+.\end{equation}
Let us go back to the $sl(n+1)$-module $\wht M$ and its submodule $\ol{M}$ defined in (1.4)-(1.11). According to (1.6)-(1.9),
\begin{equation}\msr G_+(1\otimes M)=\{0\},\qquad U(\msr G_0)(1\otimes M)= 1\otimes M. \end{equation}
By the PBW theorem,
\begin{equation}U(sl(n+1))(1\otimes M)=U(\msr G_-)U(\msr G_0)U(\msr G_+)(1\otimes M)=U(\msr G_-)(1\otimes M).\end{equation}
Furthermore, (1.11) and the first equation in (1.8) imply
\begin{equation}\ol{M}_r=\msr G_-^r(1\otimes M).\end{equation}

First we consider the case $n=2$. For $k\in\mbb{N}$, the $k$-dimensional  $sl(2)$-module
\begin{equation}
	M=\sum\limits_{i=0}^k\mathbb{C}y^i\subset\mathbb{C}[y],
\end{equation}
with the representation given by
\begin{equation}	
	E_{1,2}|_M=y^{2}\partial_{y}-ky,~E_{2,1}|_M=-\partial_{y},~(E_{1,1}-E_{2,2})|_M=2y\partial_y-k.
\end{equation}

Now we convert $M$ into an irreducible $gl(2)$-module by letting
\begin{equation}
	(E_{1,1}+E_{2,2})|_M=c\:\mathrm{Id}_M,
\end{equation}
where $c\in\mbb C$ is an arbitrary constant. Applying (1.5)-(1.9). Note
\begin{equation}\wht M=\mbb C[x_1,x_2]\otimes_{\mbb C} M=\mbox{Span}\{x_1^{\al_1}x_2^{\al_2}y^\be\mid \al_1,\al_2\in\mbb N,\;\be\in\ol{0,k}\}\subset\mbb C[x_1,x_2,y].\end{equation}	
Then $\ol M= U(sl(3)_{-})(1\otimes M)$ is an irreducible $sl(3)$-submodule of $\wht M$ with highest weight
\begin{equation}\lambda=-(\frac{3}{2}c+\frac{k}{2})\lambda_{1}+k\lambda_{2}.\end{equation}  Moreover, it is finite-dimensional if and one if we choose $c\in\frac{-2\mathbb{N}-k}{3}$. Fix such $c$.  Denote
\begin{equation}
	k_1=-(\frac{3}{2}c+\frac{k}{2}),~k_2=k,
\end{equation} where $k_1$ can be any nonnegative integer.
We have the following full projective oscillator representation of $sl(3)$:
\begin{equation}
	E_{1,2}|_{\wht M}=x_1\partial_{x_2}+y^2\partial_y-k_2y,\;\;E_{2,1}|_{\wht M}=x_2\partial_{x_1}-\partial_y,\end{equation}
\begin{equation}E_{1,3}|_{\wht M}=x_1(x_1\partial_{x_1}+x_2\partial_{x_2}-k_1)+(x_1+x_2y)(y\partial_y-k_2),\end{equation}
\begin{equation}E_{2,3}|_{\wht M}=x_2(x_1\partial_{x_1}+x_2\partial_{x_2}-k_1)-(x_1+x_2y)\partial_y,\end{equation}
\begin{equation}E_{3,1}|_{\wht M}=-\partial_{x_1},\qquad E_{3,2}|_{\wht M}=-\partial_{x_2},\end{equation}
\begin{equation}(E_{1,1}-E_{2,2})|_{\wht M}=x_1\partial_{x_1}-x_2\partial_{x_2}+2y\partial_y-k_2,\end{equation}
\begin{equation}(E_{2,2}-E_{3,3})|_{\wht M}=x_1\partial_{x_1}+2x_2\partial_{x_2}-y\partial_y-k_1.	\end{equation}	
Recall (1.11) and (1.12). We have:\pse

{\bf Proposition 2.1}\quad {\it For $0\leq r\leq k_1$, $\ol{M}_r=\wht M_r$. When $k_1<r\leq k_1+k_2$,
	\begin{equation}{\cal B}_r=\left\lbrace x_{1}^{\alpha_{1}}x_{2}^{\alpha_{2}}(x_{1}+x_{2}y)^{r-k_1}y^{\beta} | \al_1,\al_2\in\mbb N,\; \alpha_{1}+\alpha_{2}=k_1; \beta\in\ol{0,k_1+k_2-r} \right\rbrace	\end{equation}
	forms a basis of $\ol{M}_r$. If $r>k_1+k_2$, $\ol{M}_r=\{0\}$. In particular,
	\begin{equation}\left\lbrace x_{1}^{\alpha_{1}}x_{2}^{\alpha_{2}}y^{\beta} |\al_1,\al_2\in\mbb N,\; \alpha_{1}+\alpha_{2}\leq k_1; \beta\in\ol{0, k_2} \right\rbrace\bigcup(\bigcup_{r=k_1+1}^{k_1+k_2}{\cal B}_r)\end{equation}
	is a basis of the irreducible $sl(3)$-module $V_2(k_1\lmd_1+k_2\lmd_2)=\ol M$.
}
\begin{proof}
	We will prove the theorem by induction on $r$. If $r=0$, $\ol M_0=1\otimes M=\wht M_0$.
	Suppose $\ol{M}_r=\wht M_r$ for $r<\ell \leq k_1$. Consider $r=\ell$.
	For any $(\alpha_1,\alpha_2,\beta)\in\mathbb{N}^3$ with $\alpha_1+\alpha_2=\ell-1$ and $\beta\leqslant k_2-1$, we have \begin{equation}x_1^{\alpha_1}x_2^{\alpha_2}y^{\beta},x_1^{\alpha_1}x_2^{\alpha_2}y^{\beta+1}\in\ol{M}_{\ell-1}.\end{equation}
	Moreover, (2.21) and (2.22) give			
	\begin{eqnarray}& &					\begin{pmatrix}E_{13}(x^{\alpha_{1}}_{1}x^{\alpha_{2}}_{2}y^{\beta}) \\ E_{23}(x^{\alpha_{1}}_{1}x^{\alpha_{2}}_{2}y^{\beta+1})
		\end{pmatrix}
		\nonumber \\ &=&
		\begin{pmatrix}
			(\ell+\beta-k_1-k_2-1) & \beta-k_2 \\ -\beta-1 & \ell-\beta-k_1-2
		\end{pmatrix}
		\begin{pmatrix}
			x_{1}^{\alpha_{1}+1}x_{2}^{\alpha_{2}}y^{\beta}\\x_{1}^{\alpha_{1}}x_{2}^{\alpha_{2}+1}y^{\beta+1}
		\end{pmatrix}
	\end{eqnarray}
	and
	\begin{equation}		
		x_{1}^{\alpha_{1}+1}x_{2}^{\alpha_{2}}y^{k_2}=\frac{1}{\ell-k_1-1}E_{13}(x^{\alpha_{1}}_{1}x^{\alpha_{2}}_{2}y^{k_2})\in\ol{M}_\ell,\end{equation}
	\begin{equation}					x_{1}^{\alpha_{1}}x_{2}^{\alpha_{2}+1}=\frac{1}{\ell-k_1-1}E_{23}(x^{\alpha_{1}}_{1}x^{\alpha_{2}}_{2})\in\ol{M}_\ell.				\end{equation}
	Since
	\begin{equation}
		\begin{vmatrix}
			(\ell+\beta-k_1-k_2-1) & \beta-k_2 \\ -\beta-1 & \ell-\beta-k_1-2
		\end{vmatrix}
		=(\ell-k_1-1)(\ell-k_1-k_2-2)\neq0,
	\end{equation}
	Solving (2.29) yields	
	\begin{equation}x_{1}^{\alpha_{1}+1}x_{2}^{\alpha_{2}}y^{\beta},x_{1}^{\alpha_{1}}x_{2}^{\alpha_{2}+1}y^{\beta+1}\in \ol{M}_\ell.\end{equation}
	According to (2.30)-(2.33), $\ol{M}_\ell=\wht M_\ell$. By induction, $\ol{M}_r=\wht M_r$ for $r\in\ol{0,k_1}$.
	
	For technical convenience, we allow $r=k_1$ in (2.26). Then
	\begin{equation}{\cal B}_{k_1}=\left\lbrace x_{1}^{\alpha_{1}}x_{2}^{\alpha_{2}}y^{\beta} | \al_1,\al_2\in\mbb N,\; \alpha_{1}+\alpha_{2}=k_1; \beta\in\ol{0,k_2} \right\rbrace\;\mbox{is a basis of}\;\ol{M}_{k_1}=\wht M_{k_1}.	\end{equation}
	In this case, $\msr G_-=\mbb CE_{1,3}+\mbb CE_{2,3}$ (cf. (2.9)). For an element in (2.26), we have
	\begin{eqnarray}& &E_{13}(x_1^{\alpha_1}x_{2}^{\alpha_2}(x_1+x_2y)^{r-k_1}y^{\beta})\nonumber\\
		&=&(x_1(x_1\partial_{x_1}+x_2\partial_{x_2}-k_1)+(x_1+x_2y)(y\partial_y-k_2))(x_1^{\alpha_1}x_{2}^{\alpha_2}(x_1+x_2y)^{r-k_1}y^{\beta})
		\nonumber\\ &=&(r+\be-k_1-k_2)x_1^{\alpha_1}x_{2}^{\alpha_2}(x_1+x_2y)^{r+1-k_1}y^{\beta} \end{eqnarray}
	and
	\begin{eqnarray}& &E_{23}(x_1^{\alpha_1}x_{2}^{\alpha_2}(x_1+x_2y)^{r-k_1}y^{\beta})\nonumber\\
		&=&(x_2(x_1\partial_{x_1}+x_2\partial_{x_2}-k_1)-(x_1+x_2y)\partial_y)
		(x_1^{\alpha_1}x_{2}^{\alpha_2}(x_1+x_2y)^{r-k_1}y^{\beta})
		\nonumber\\ &=&-\be x_1^{\alpha_1}x_{2}^{\alpha_2}(x_1+x_2y)^{r+1-k_1}y^{\beta} \end{eqnarray}
	by (2.21) and (2.22). Note that if $\be=k_1+k_2-r$, then (2.35) becomes 0. Thus if ${\cal B}_r$ is a basis of $\ol{M}_r$,
	(2.35) and (2.36) implies that ${\cal B}_{r+1}$ is a basis of $\ol{M}_{r+1}$ by (2.13). Starting from $r=k_1$ and (2.34), we obtain that
	${\cal B}_r$ is a basis of $\ol{M}_r$ for any $r\in\ol{k_1+1,k_1+k_2}$ by induction. Since
	\begin{equation}{\cal B}_{k_1+k_2}=\left\lbrace x_{1}^{\alpha_{1}}x_{2}^{\alpha_{2}}(x_{1}+x_{2}y)^{k_2} | \al_1,\al_2\in\mbb N,\; \alpha_{1}+\alpha_{2}=k_1 \right\rbrace,	\end{equation}		
	(2.35) and (2.36) with $r=k_1+k_2$ and $\be=0$ imply 		
	\begin{equation}E_{13}({\cal  B}_{k_1+k_2})=E_{23}({\cal B}_{k_1+k_2})={0}. \end{equation}		
	By (2.13), $\ol M_r=\{0\}$ for $k_1+k_2<r\in\mbb N$.
	
\end{proof}

Note
\begin{equation}
	\dim \ol{M}_{r}=\left\lbrace
	\begin{array}{ll}
		(r+1)(k_2+1) &\mathrm{if}~0\leqslant r \leqslant k_1, \\
		(k_1+1)(k_1+k_2+1-r) &\mathrm{if}~k_1<r \leqslant k_1+k_2.
	\end{array}
	\right.
\end{equation}
Moreover, (1.17) gives
\begin{equation}
	d_2(k_1\lmd_1+k_2\lmd_2)=\frac{(k_1+1)(k_2+1)(k_1+k_2+2)}{2}.
\end{equation}
Thus the equation $\sum_{r=0}^{k_1+k_2}\dim\ol M_r=\dim\ol M=\dim V_2(\lmd)$
becomes
\begin{equation}\sum_{r=0}^{k_1}(r+1)(k_2+1)+\sum_{r=k_1+1}^{k_1+k_2}(k_1+1)(k_1+k_2+1-r)=\frac{(k_1+1)(k_2+1)(k_1+k_2+2)}{2}.
\end{equation}
This completes the realization of the irreducible $sl(3)$-module $V_2(k_1\lmd_1+k_2\lmd_2)=\ol{M}$ in the variable $y,x_1,x_2$ (cf. (2.18) and (2.19)).\psp

Let
\begin{equation}\msr C_n=\mathbb{C}[x_{i,j}|1\leqslant j \leqslant i\leqslant n]\end{equation}	
be the  polynomial algebra in $n(n+1)/2$ variables. For the convenience, we make the following convention	

\begin{equation}
	x_{i,j}=\left\lbrace
	\begin{aligned}
		&\mathrm{variable}~~x_{i,j},   ~&1\leqslant j \leqslant i, \\
		&-1, ~&j=i+1,\\
		&0, ~&\mbox{otherwise}.
	\end{aligned}
	~~~~~~~~~i\in\overline{1,n},
	\right.
\end{equation}
Moreover, we take the notations
\begin{equation}
	\alpha_i=(\alpha_{i,1},\alpha_{i,2},...,\alpha_{i,i})\in \mathbb{N}^{i},~~i\in\overline{1,n}
\end{equation}
and
\begin{equation}
	\alpha=(\alpha_{1},\alpha_{2},..,\alpha_{n})\in\mathbb{N}^{\frac{n(n+1)}{2}}
\end{equation}
Denote by $X_{i}^{\alpha_{i}}$ and $X^{\alpha}$ the monomials
\begin{equation}
	X_i^{\alpha_{i}}=x_{i,1}^{\alpha_{i,1}}x_{i,2}^{\alpha_{i,2}}...x_{i,i}^{\alpha_{i,i}},
\end{equation}
and
\begin{equation}
	X^{\alpha}=X_1^{\alpha_1}X_2^{\alpha_2}...X_n^{\alpha_n}.
\end{equation}
\psp

For $1\leqslant j\leqslant i\leqslant n$ and $\Theta=(\theta_1,\theta_2,...,\theta_{i-j+1})\in \mathbb{N}^{i-j+1}$, we denote
\begin{equation}
	d_{i,j}(\Theta)=\begin{vmatrix}
		x_{j,\theta_{1}}&x_{j,\theta_{2}}&...&x_{j,\theta_{i-j+1}}\\
		x_{j+1,\theta_{1}}&x_{j+1,\theta_{2}}&...&x_{j+1,\theta_{i-j+1}}\\
		...&...&...&...\\
		x_{i,\theta_{1}}&x_{i,\theta_{2}}&...&x_{i,\theta_{i-j+1}}\\
	\end{vmatrix}
\end{equation}
In particular, we take
\begin{equation}
	\begin{aligned}
		D_{i,j}(s)=&d_{i,j}(s,j+1,j+2,...,i)\pse\\
		=&\begin{vmatrix}
			x_{j,s} &-1&0&...&0 \\
			x_{j+1,s}&x_{j+1,j+1}&-1&...&0 \\
			...&...&...&...&...\\
			x_{i-1,s}&x_{i-1,j+1}&x_{i-1,j+2}&...&-1\\
			x_{i,s}&x_{i,j+1}&x_{i,j+2}&...&x_{i,i}
		\end{vmatrix}.
	\end{aligned}
\end{equation}\psp

{\bf Lemma 2.2}\quad {\it
	For $1\leqslant s\leqslant j<i\leqslant n$, the following equations hold,
	
	(i) $D_{i,j}(s)=x_{i,s}+\sum\limits_{t=j+1}^ix_{i,t}D_{t-1,j}(s)$;
	
	(ii) $D_{i,j}(s)=x_{i,s}+\sum\limits_{t=j}^{i-1}x_{t,s}D_{i,t+1}(t+1)$.
}
\begin{proof}
	
	The algebraic cofactor of $x_{i,t},(j+1\leqslant t\leqslant i)$ in $D_{i,j}(s)$ is
	\begin{eqnarray}
		A_{i,t}&=&(-1)^{i+t}\times\nonumber\\ & &\begin{vmatrix}
			x_{j,s}&-1&0&...&0&0&0&...&0\\
			x_{j+1,s}&x_{j+1,j+1}&-1&...&0&0&0&...&0\\
			...&...&...&...&...&...&...&...&...\\
			x_{t-1,s}&x_{t-1,j+1}&x_{t-1,j+2}&...&x_{t-1,t-1}&0&0&...&0\\
			x_{t,s}&x_{t,j+1}&x_{t,j+2}&...&x_{t,t-1}&-1&0&...&0\\
			x_{t+1,s}&x_{t+1,j+1}&x_{t+1,j+2}&...&x_{t+1,t-1}&x_{t+1,t+1}&-1&...&0\\
			...&...&...&...&...&...&...&...&...\\
			x_{i-1,s}&x_{i-1,j+1}&x_{i-1,j+2}&...&x_{i-1,t-1}&x_{i-1,t+1}&x_{i-1,t+2}&...&-1
		\end{vmatrix}	\nonumber \\
		&=&\begin{vmatrix}
			x_{j,s}&-1&0&...&0\\
			x_{j+1,s}&x_{j+1,j+1}&-1&...&0\\
			...&...&...&...&...\\
			x_{t-1,s}&x_{t-1,j+1}&x_{t-1,j+2}&...&x_{t-1,t-1}
		\end{vmatrix}\nonumber \\
		&=&D_{t-1,j}(s).
	\end{eqnarray}
	Hence expanding the determinant (2.49) according to the last row, we get
	\begin{equation}D_{i,j}(s)=x_{i,s}+\sum\limits_{t=j+1}^ix_{i,t}A_{i,t}=x_{i,s}
		+\sum\limits_{t=j+1}^ix_{i,t}D_{t-1,j}(s).\end{equation}
	which proves (i). Expanding the determinant (2.49) according to the first  row, we obtain
	\begin{eqnarray}D_{i,j}(s)&=&x_{j,s}D_{i,j+1}(j+1)+D_{i,j+1}(s)\nonumber\\
		&=&x_{j,s}D_{i,j+1}(j+1)+x_{j+1,s}D_{i,j+2}(j+2)+D_{i,j+2}(s)\nonumber\\ &=&\cdots
		=\sum\limits_{t=j}^{i-1}x_{t,s}D_{i,t+1}(t+1)+D_{i,i}(s).\end{eqnarray}
	Note that $D_{i,i}(s)=x_{i,s}$, Conclusion (ii) holds.
\end{proof}

Recall our realization of the irreducible $sl(3)$-module $V_2(k_1\lmd_1+k_2\lmd_2)$ in (2.14)-(2.41).
We redenote
\begin{equation}x_{1,1}=y,\;\;x_{2,1}=x_1,\;\;x_{2,2}=x_2.\end{equation} Thus we have used (1.4)-(1.12) to realized the
$sl(3)$-module $V_2(k_1\lmd_1+k_2\lmd_2)$ as a $sl(3)$-submodule of the  $sl(3)$-module $\msr C_2=\mbb C[x_{1,1},x_{2,1},x_{2,2}]$, whose representation formulas are given in (2.20)-(2.25) with $\wht M$ replaced by $\msr C_2$. Suppose that we have used (1.4)-(1.12) successively to realize the irreducible $sl(n)$-module $M=V_{n-1}(k_2\lambda_1+k_3\lambda_2+\cdots+k_n\lambda_{n-1})$ as a submodule of the $sl(n)$-module $\msr C_{n-1}=\mbb C[x_{i,j}|1\leqslant j \leqslant i\leqslant n-1]$.
Fix $c_n\in\mbb C$ and impose
\begin{equation}I_n|_M=c_n\mathrm{Id}_M.\end{equation}	
Then $M$ become an irreducible $gl(n)$-module. Use (1.4)-(1.12) with
\begin{equation}x_1=x_{n,1},\;\;x_2=x_{n,2},\;\;...,\;\;x_n=x_{n,n}.\end{equation}
Now
\begin{eqnarray}\wht M&=&\mbb C[x_{n,1},x_{n,2},...,x_{n,n}]\otimes_{\mbb C} M\nonumber\\ &\subset&
	\mbb C[x_{n,1},x_{n,2},...,x_{n,n}]\otimes_{\mbb C}\msr C_{n-1}\nonumber\\
	&=&\mbb C[x_{n,1},x_{n,2},...,x_{n,n}]\otimes_{\mbb C}
	\mbb C[x_{i,j}|1\leqslant j \leqslant i\leqslant n-1]\nonumber\\
	&=&\mbb C[x_{i,j}|1\leqslant j \leqslant i\leqslant n]=\msr C_n.\end{eqnarray}
Then $\ol M=V_n(\lmd)$
is an irreducible module with highest weight \begin{equation}\lambda=-\frac{1}{n}((n+1)c_n+\sum\limits_{i=1}^{n-1}(n-i)k_{i+1})\lambda_1
	+k_2\lambda_2+k_3\lambda_3+\cdots+k_n\lambda_n.\end{equation}	 In particular, we can choose appropriate $c_n\in \mbb C$ such that \begin{equation}k_1=-\frac{1}{n}((n+1)c_n+\sum\limits_{i=1}^{n-1}(n-i)k_{i+1})\in\mathbb{N}.\end{equation}	
Thus $\ol M$ is a finite-dimensional irreducible $sl(n+1)$-module. Indeed, $k_1$ can take any nonnegative integer.

For convenience, we extend the representation of $sl(n+1)$ on $\wht M$ to the representation of $sl(n+1)$ on $\msr C_n$ by (1.4)-(1.9) with $M$ replaced by $\msr C_{n-1}$ and $\wht{M}$ by $\msr C_n$.
Fix $c_{n+1}\in\mbb C$. We make $\msr C_n$ as a $gl(n+1)$-module by imposing
\begin{equation}I_{n+1}|_{\msr C_n}=c_{n+1}\mathrm{Id}_{\msr C_n}.\end{equation}

{\bf Theorem 2.3}\quad{\it
	Take the convention in (2.43). The representation formulas of $gl(n+1)$ on $\msr C_n$ are as follows:
	\begin{equation}E_{i,i}|_{\msr C_n}=\frac{1}{n+1}c_{n+1}-c_n+\sum\limits_{r=i}^{n}(x_{r,i}\partial_{x_{r,i}}-k_{n-r+1})-\sum\limits_{s=1}^{i-1}x_{i-1,s}\partial_{x_{i-1,s}}	\;\;\for\;\;i\in\ol{1,n+1},\end{equation}
	\begin{equation}E_{i,j}|_{\msr C_n}=\sum\limits_{r=i-1}^{n}x_{r,i}\partial_{x_{r,j}} \qquad\for\;\;1\leqslant j<i\leqslant n+1 \end{equation}
	and for $1\leqslant i<j\leqslant n+1$,
	\begin{equation}
		E_{i,j}|_{\msr C_n}=\sum\limits_{r=j}^{n}x_{r,i}\partial_{x_{r,j}}+\sum\limits_{s=i-1}^{j-1}\sum\limits_{t=1}^{s}D_{j-1,s}(t)x_{s,i}\partial_{x_{s,t}}-\sum\limits_{s=i}^{j-1}D_{j-1,s}(i)k_{n+1-s}.	\end{equation}}
\begin{proof}
	When $n=2$, (2.60)-(2.62) become
	\begin{equation}E_{12}|_{\mathscr{C}_2}=x_{2,1}\partial_{x_{2,2}}+x_{1,1}^2\partial_{x_{1,1}}-k_2x_{1,1},\;\;E_{21}|_{\mathscr{C}_2}=x_{2,2}\partial_{x_{2,1}}-\partial_{x_{1,1}},\end{equation}	\begin{equation}E_{13}|_{\mathscr{C}_2}=x_{2,1}(x_{2,1}\partial_{x_{2,1}}+x_{2,2}\partial_{x_{2,2}}-k_1)+
		(x_{2,1}+x_{2,2}x_{1,1})(x_{1,1}\partial_{x_{1,1}}-k_2),\end{equation}
	\begin{equation}		E_{23}|_{\mathscr{C}_2}=x_{2,2}(x_{2,1}\partial_{x_{2,1}}+x_{2,1}\partial_{x_{2,1}}-k_1)-(x_{2,1}+x_{2,2}x_{1,1})
		\partial_{x_{1,1}},\end{equation}
	\begin{equation}E_{31}|_{\mathscr{C}_2}=-\partial_{x_{2,1}},~E_{32}|_{\mathscr{C}_2}=-\partial_{x_{2,2}},\end{equation}
	\begin{equation}E_{11}|_{\mathscr{C}_2}=x_{2,1}\partial_{x_{2,1}}+x_{1,1}\partial_{x_{1,1}}+\frac{1}{3}c_3-c_2-k_1-k_2,\end{equation} \begin{equation}E_{22}|_{\mathscr{C}_2}=x_{2,2}\partial_{x_{2,2}}-x_{1,1}\partial_{x_{1,1}}+\frac{1}{3}c_3-c_2-k_1,\end{equation}
	\begin{equation}E_{33}|_{\mathscr{C}_2}=-x_{2,1}\partial_{x_{2,1}}-x_{2,2}\partial_{x_{2,2}}+\frac{1}{3}c_3-c_2.\end{equation}	Under the identification (2.53), (2.63)-(2.66) are exactly (2.20)-(2.23). Moreover, (2.67)$-$(2.68) is (2.24) and (2.68)$-$(2.69) is (2.25). According to (2.19), the sum of (2.67)-(2.69) is exactly (2.59) with $n=2$.	So the theorem holds for the case of $gl(3)$.

	Suppose that the theorem holds for $sl(n)$.  In the rest of this paper, we make a convention that
	\begin{equation}\mbox{differential operators on}\; {\msr C}_{n-1}\;\mbox{are also viewed as those on}\;\msr C_n.\end{equation}
	Let $i\in\ol{1,n}$. First
	\begin{equation}
		(E_{i,i}-E_{n+1,n+1})|_{\msr C_n}= \sum\limits_{r=1}^nx_{n,r}\partial_{x_{n,r}}+x_{n,i}\partial_{x_{n,i}}+
		(I_n+E_{i,i})|_{\msr C_{n-1}}
	\end{equation}
	by (1.7) and (1.9). Note $M=V_{n-1}(k_2\lmd_1+k_3\lmd_2+\cdots+k_n\lmd_{n-1})$.  Applying (2.60) with $n$ replaced by $n-1$,  we get
	\begin{equation}
		(I_n+E_{i,i})|_{\msr C_{n-1}}=c_n+\frac{1}{n}c_n-c_{n-1}+\sum\limits_{s=i}^{n-1}(x_{s,i}\partial_{x_{s,i}}-k_{n-s+1})-\sum\limits_{s=1}^{i-1}x_{i-1,s}\partial_{x_{i-1,s}}.
	\end{equation}
	
	According to (2.58),
	\begin{equation}c_n=-\frac{1}{n+1}(nk_1+\sum\limits_{i=1}^{n-1}(n-i)k_{i+1}).\end{equation}	
	By our inductive construction,
	\begin{equation}c_{n-1}=-\frac{1}{n}((n-1)k_2+\sum\limits_{i=1}^{n-2}(n-1-i)k_{i+2}).\end{equation}	
	Thus
	\begin{equation}\frac{n+1}{n}c_n-c_{n-1}=-k_1.\end{equation}	
	Expressions (2.71)-(2.72) and (2.75) imply		
	\begin{equation}
		(E_{i,i}-E_{n+1,n+1})_{\msr C_n}=\sum\limits_{r=1}^{n}x_{n,r}\partial_{x_{n,r}}+\sum\limits_{s=i}^{n}(x_{s,i}\partial_{x_{s,i}}-k_{n-s+1})
		-\sum\limits_{s=1}^{i-1}x_{i-1,s}\partial_{x_{i-1,s}}.\end{equation}
	Hence	
	\begin{eqnarray}		
		E_{n+1,n+1}|_{\msr C_n}&=&\frac{1}{n+1}( \sum\limits_{r=1}^{n+1}E_{r,r}-\sum\limits_{r=1}^{n}(E_{r,r}-E_{n+1,n+1}))|_{\msr C_n}\nonumber\\
		&=&\frac{1}{n+1}c_{n+1}-c_n-\sum\limits_{r=1}^{n}x_{n,r}\partial_{x_{n,r}}
	\end{eqnarray}
	by (2.76), and		
	\begin{eqnarray}E_{i,i}|_{\msr C_n}&=&(E_{i,i}-E_{n+1,n+1}+E_{n+1,n+1})|_{\msr C_n}\nonumber\\					&=&\frac{1}{n+1}c_{n+1}-c_n+\sum\limits_{s=i}^{n}(x_{s,i}\partial_{x_{s,i}}-k_{n-s+1})
		-\sum\limits_{s=1}^{i-1}x_{i-1,s}\partial_{x_{i-1,s}}.	\end{eqnarray}

	Assume $1\leqslant j<i\leqslant n+1$.		
	If $i=n+1$,
	\begin{equation}
		E_{n+1,j}|_{\msr C_n}=-\partial_{x_{n,j}}=x_{n,n+1}\partial_{x_{n,j}},
	\end{equation}
	where the last equation dues to the convention (2.43) that $x_{n,n+1}=-1$. When $i<n+1$,
	\begin{equation}
		E_{i,j}|_{\msr C_n}=x_{n,i}\partial_{x_{n,j}}+E_{i,j}|_{\msr C_{n-1}}=\sum\limits_{r=i-1}^nx_{r,i}\partial_{x_{r,j}},
	\end{equation}
	where the last equation is from inductive assumption.
	
	Let $1\leqslant i<j\leqslant n+1$.	We want to prove (2.62). Note that
	\begin{equation}\sum\limits_{r=i+1}^n\sum\limits_{s=i}^{r-1}\sum\limits_{t=1}^{s}
		=\sum\limits_{s=i}^{n-1}\sum\limits_{t=1}^{s}\sum\limits_{r=s+1}^{n},\qquad \sum\limits_{r=i+1}^n\sum\limits_{s=i}^{r-1}
		=\sum\limits_{s=i}^{n-1}\sum\limits_{r=s+1}^{n}. \end{equation}
	
	Consider $j=n+1$. By inductive assumption and (2.81), we have
	\begin{eqnarray}& &
		E_{i,n+1}|_{\msr C_n}\nonumber\\&=&x_{n,i}\sum\limits_{r=1}^nx_{n,r}\partial_{x_{n,r}}+ x_{n,i}I_n|_{\msr C_{n-1}}+\sum\limits_{r=1}^nx_{n,r} E_{i,r}|_{\msr C_{n-1}}\nonumber\\&=&x_{n,i}(c_n+\sum\limits_{r=1}^nx_{n,r}\partial_{x_{n,r}})+\sum\limits_{r=1}^nx_{n,r} E_{i,r}|_{\msr C_{n-1}}\nonumber\\ &=&x_{n,i}(c_n+\sum\limits_{r=1}^nx_{n,r}\partial_{x_{n,r}})+
		\sum\limits_{r=1}^{i-1}x_{n,r}E_{i,r}|_{\msr C_{n-1}}+x_{n,i}E_{i,i}|_{\msr C_{n-1}}+\sum\limits_{r=i+1}^{n}x_{n,r}E_{i,r}|_{\msr C_{n-1}}
		\nonumber\\&=&x_{n,i}(c_n+\sum\limits_{r=1}^nx_{n,r}\partial_{x_{n,r}})+
		\sum\limits_{r=1}^{i-1}\sum\limits_{s=i-1}^{n-1}x_{n,r}x_{s,i}\partial_{x_{s,r}}\nonumber\\
		&&+x_{n,i}\left(\frac{1}{n}c_n-c_{n-1}+\sum\limits_{s=i}^{n-1}(x_{s,i}\partial_{x_{s,i}}-k_{n-s+1})
		-\sum\limits_{s=1}^{i-1}x_{i-1,s}\partial_{x_{i-1,s}}\right)\nonumber\\		&&+\sum\limits_{r=i+1}^{n}x_{n,r}\left(\sum\limits_{s=r}^{n-1}x_{s,i}\partial_{x_{s,r}}
		+\sum\limits_{s=i-1}^{r-1}\sum\limits_{t=1}^{s}D_{r-1,s}(t)x_{s,i}\partial_{x_{s,t}}
		-\sum\limits_{s=i}^{r-1}D_{r-1,s}(i)k_{n-s+1}\right)
		\nonumber\\				&=&\!\!\!\!\sum\limits_{s=i}^{n-1}\sum\limits_{t=1}^{s}\left(x_{n,t}+\sum\limits_{r=s+1}^{n}x_{n,r}D_{r-1,s}(t)\right)
		x_{s,i}\partial_{x_{s,t}}-\sum\limits_{s=i}^{n-1}\left(x_{n,i}+\sum\limits_{r=s+1}^{n}x_{n,r}D_{r-1,s}(i)\right)k_{n-s+1}\nonumber\\
		&&-\sum\limits_{t=1}^{i-1}\left(x_{n,t}+x_{n,i}x_{i-1,t}
		+\sum\limits_{r=i+1}^{n}x_{n,r}D_{r-1,i-1}(t)\right)\partial_{x_{i-1,t}}
		-k_1x_{n,i}+x_{n,i}\sum\limits_{r=1}^nx_{n,r}\partial_{x_{n,r}}
		\nonumber\\				&=&\!\!\!\!\sum\limits_{s=i-1}^n\sum\limits_{t=1}^{s}\left(x_{n,t}+\sum\limits_{r=s+1}^{n}x_{n,r}D_{r-1,s}(t)\right)
		x_{s,i}\partial_{x_{s,t}}-\sum\limits_{s=i}^n\left(x_{n,i}+\sum\limits_{r=s+1}^{n}x_{n,r}D_{r-1,s}(i)\right)k_{n-s+1}
		\nonumber\\	
		&=&\sum\limits_{s=i-1}^n\sum\limits_{t=1}^{s}D_{n,s}(t)x_{s,i}\partial_{x_{s,t}}-\sum\limits_{s=i}^{n}D_{n,s}(i)k_{n-s+1},
	\end{eqnarray}
	where in the last equality, we have used Lemma 2.1 (i). When $j<n+1$, inductive assumption and (2.62) with $n$ replaced by $n-1$ imply
	\begin{eqnarray}
		E_{i,j}|_{\msr C_n}&=&x_{n,i}\partial_{x_{n,j}}+ E_{i,j}|_{\msr C_{n-1}}\nonumber\\
		&=&x_{n,i}\partial_{x_{n,j}}+\sum\limits_{r=j}^{n-1}x_{r,i}\partial_{x_{r,j}}+
		\sum\limits_{s=i-1}^{j-1}\sum\limits_{t=1}^{s}D_{j-1,s}(t)x_{s,i}\partial_{x_{s,t}}
		-\sum\limits_{s=i}^{j-1}D_{j-1,s}(i)k_{n+1-s} \nonumber\\&=&
		\sum\limits_{r=j}^{n}x_{r,i}\partial_{x_{r,j}}+\sum\limits_{s=i-1}^{j-1}\sum\limits_{t=1}^{s}D_{j-1,s}(t)
		x_{s,i}\partial_{x_{s,t}}-\sum\limits_{s=i}^{j-1}D_{j-1,s}(i)k_{n+1-s}.\end{eqnarray}
\end{proof}

Next we will find the presentation of the $sl(n+1)$-module $\ol M=V(\lmd)$ (cf. (1.10)-(1.12) and (2.55)).
Recall again that  $x_{s,i}=0$ if $s<i-1$ by the convention (2.43), and $D_{n,s}(i)=0$ when $s<i\leq n$ by (2.49). Thus by the last equality in (2.82),  we can write
\begin{equation}
	E_{i,n+1}|_{\msr C_n}=\Phi_i-\Psi_i\;\;\mbox{with}\;\;\Phi_i=\sum\limits_{1\leqslant t\leqslant s\leqslant n}D_{n,s}(t)x_{s,i}\partial_{x_{s,t}},\;\;\Psi_i=\sum\limits_{s=1}^{n}D_{n,s}(i)k_{n-s+1}. \label{3.26}
\end{equation}
Now we calculate
\begin{eqnarray}
	& &\Phi_i(d_{s,t}(\Theta))=\Phi_i\left( \begin{vmatrix}
		x_{t,\theta_1}&x_{t,\theta_2}&...&x_{t,\theta_{s-t+1}}\\
		x_{t+1,\theta_1}&x_{t+1,\theta_2}&...&x_{t+1,\theta_{s-t+1}}\\
		...&...&...&...\\
		x_{s,\theta_1}&x_{s,\theta_2}&...&x_{s,\theta_{s-t+1}}
	\end{vmatrix}\right) \nonumber\\
	&=&\sum\limits_{j=t}^sx_{j,i}\begin{vmatrix}
		x_{t,\theta_1}&x_{t,\theta_2}&...&x_{t,\theta_{s-t+1}}\\
		...&...&...&...\\
		x_{j-1,\theta_1}&x_{j-1,\theta_2}&...&x_{j-1,\theta_{s-t+1}}\\
		D_{n,j}(\theta_1)&D_{n,j}(\theta_2)&..&D_{n,j}(\theta_{s-t+1})\\
		x_{j+1,\theta_1}&x_{j+1,\theta_2}&...&x_{j+1,\theta_{s-t+1}}\\
		...&...&...&...\\
		x_{s,\theta_1}&x_{s,\theta_2}&...&x_{s,\theta_{s-t+1}}
	\end{vmatrix} \nonumber\\
	&=&\sum\limits_{j=t}^sx_{j,i}\begin{vmatrix}
		x_{t,\theta_1}&x_{t,\theta_2}&...&x_{t,\theta_{s-t+1}}\\
		...&...&...&...\\
		x_{j-1,\theta_1}&x_{j-1,\theta_2}&...&x_{j-1,\theta_{s-t+1}}\\		\sum\limits_{r=j}^{n}x_{r,\theta_1}D_{n,r+1}(r+1)&\sum\limits_{r=j}^{n}x_{r,\theta_2}D_{n,r+1}(r+1)&..&\sum\limits_{r=j}^{n}x_{r,\theta_{s-t+1}}D_{n,r+1}(r+1)\\
		x_{j+1,\theta_1}&x_{j+1,\theta_2}&...&x_{j+1,\theta_{s-t+1}}\\
		...&...&...&...\\
		x_{s,\theta_1}&x_{s,\theta_2}&...&x_{s,\theta_{s-t+1}}
	\end{vmatrix}\nonumber\\
	&=&\sum\limits_{j=t}^sx_{j,i}\begin{vmatrix}
		x_{t,\theta_1}&x_{t,\theta_2}&...&x_{t,\theta_{s-t+1}}\\
		...&...&...&...\\
		x_{j-1,\theta_1}&x_{j-1,\theta_2}&...&x_{j-1,\theta_{s-t+1}}\\
		x_{j,\theta_1}D_{n,j+1}(j+1)&x_{j,\theta_2}D_{n,j+1}(j+1)&..&x_{j,\theta_{s-t+1}}D_{n,j+1}(j+1)\\
		x_{j+1,\theta_1}&x_{j+1,\theta_2}&...&x_{j+1,\theta_{s-t+1}}\\
		...&...&...&...\\
		x_{s,\theta_1}&x_{s,\theta_2}&...&x_{s,\theta_{s-t+1}}
	\end{vmatrix}\nonumber\end{eqnarray}\begin{eqnarray}
	& +&\!\!\!\!\!\sum\limits_{j=t}^sx_{j,i}\begin{vmatrix}
		x_{t,\theta_1}&x_{t,\theta_2}&...&x_{t,\theta_{s-t+1}}\\
		...&...&...&...\\
		x_{j-1,\theta_1}&x_{j-1,\theta_2}&...&x_{j-1,\theta_{s-t+1}}\\		\sum\limits_{r=s+1}^{n}x_{r,\theta_1}D_{n,r+1}(r+1)&\sum\limits_{r=s+1}^{n}x_{r,\theta_2}D_{n,r+1}(r+1)&..&\sum\limits_{r=s+1}^{n}x_{r,\theta_{s-t+1}}D_{n,r+1}(r+1)\\
		x_{j+1,\theta_1}&x_{j+1,\theta_2}&...&x_{j+1,\theta_{s-t+1}}\\
		...&...&...&...\\
		x_{s,\theta_1}&x_{s,\theta_2}&...&x_{s,\theta_{s-t+1}}
	\end{vmatrix} \nonumber\\
	\!\!\!\!\!\!\!\! \!\!\!\!\!\!\!\!  \!\!\!\!\!\!\!\! \!\!\!\!\!\!\!\! 	&=& \!\!\!\sum\limits_{j=t}^s x_{j,i}D_{n,j+1}(j+1)d_{s,t}(\Theta)\!+
	\sum\limits_{j=t}^s x_{j,i}\begin{vmatrix}
		x_{t,\theta_1}&x_{t,\theta_2}&...&x_{t,\theta_{s-t+1}}\\
		...&...&...&...\\
		x_{j-1,\theta_1}&x_{j-1,\theta_2}&...&x_{j-1,\theta_{s-t+1}}\\
		D_{n,s+1}(\theta_1)&D_{n,s+1}(\theta_2)&..&D_{n,s+1}(\theta_{s-t+1})\\
		x_{j+1,\theta_1}&x_{j+1,\theta_2}&...&x_{j+1,\theta_{s-t+1}}\\
		...&...&...&...\\
		x_{s,\theta_1}&x_{s,\theta_2}&...&x_{s,\theta_{s-t+1}}
	\end{vmatrix} \label{3.27}
\end{eqnarray}
by Lemma 2.2 (ii).  Here we make the convention $D_{n,n+1}(n+1)=1$ and $D_{n,n+1}(\theta)=0$ if $1\leqslant\theta\leqslant n$. In the following, putting  a bracket on a row in a determinant means removing that row. If $s=n$,
\begin{equation}
	(2.85)=\sum\limits_{j=t}^nx_{j,i}D_{n,j+1}(j+1)d_{n,t}(\Theta)=D_{n,t}(i)d_{n,t}(\Theta). \label{3.28}
\end{equation}
Suppose $s=n-1$. Then we have
\begin{eqnarray}
	(2.85)&=&\sum\limits_{j=t}^{n-1}x_{j,i}D_{n,j+1}(j+1)d_{n-1,t}(\Theta)+\sum\limits_{j=t}^{n-1}x_{j,i}\begin{vmatrix}
		x_{t,\theta_1}&x_{t,\theta_2}&...&x_{t,\theta_{n-t}}\\
		...&...&...&...\\
		x_{j-1,\theta_1}&x_{j-1,\theta_2}&...&x_{j-1,\theta_{n-t}}\\
		x_{n,\theta_1}&x_{n,\theta_2}&..&x_{n,\theta_{n-t}}\\
		x_{j+1,\theta_1}&x_{j+1,\theta_2}&...&x_{j+1,\theta_{n-t+}}\\
		...&...&...&...\\
		x_{n-1,\theta_1}&x_{n-1,\theta_2}&...&x_{n-1,\theta_{n-t}}
	\end{vmatrix}\nonumber\\
	&=&D_{n,t}(i)d_{s,t}(\Theta)-x_{n,i}d_{s,t}(\Theta)+\sum\limits_{j=t}^s(-1)^{n-j-1}x_{j,i}\begin{vmatrix}
		x_{t,\theta_1}&x_{t,\theta_2}&...&x_{t,\theta_{n-t}}\\
		...&...&...&...\\
		(x_{j,\theta_1}&x_{j,\theta_2}&...&x_{j,\theta_{n-t}})\\
		...&...&...&...\\
		x_{n,\theta_1}&x_{n,\theta_2}&...&x_{n,\theta_{n-t}}\\
	\end{vmatrix}\nonumber\\
	&=&D_{n,t}(i)d_{s,t}(\Theta)+(-1)^{n-t-1}d_{n,t}(i,\theta_1,...,\theta_{n-t}).
	\label{3.29}
\end{eqnarray}

Assume $s\leqslant n-2$. Now
\begin{eqnarray}
	(2.85)&=&\sum\limits_{j=t}^sx_{j,i}D_{n,j+1}(j+1)d_{s,t}(\Theta)\nonumber\\
	& &+\sum\limits_{j=t}^{s}x_{j,i}\begin{vmatrix}
		x_{t,\theta_1}&x_{t,\theta_2}&...&x_{t,\theta_{s-t+1}}\\
		...&...&...&...\\
		x_{j-1,\theta_1}&x_{j-1,\theta_2}&...&x_{j-1,\theta_{s-t+1}}\\
		D_{n,s+1}(\theta_1)&D_{n,s+1}(\theta_2)&..&D_{n,s+1}(\theta_{s-t+1})\\
		x_{j+1,\theta_1}&x_{j+1,\theta_2}&...&x_{j+1,\theta_{s-t+1}}\\
		...&...&...&...\\
		x_{s,\theta_1}&x_{s,\theta_2}&...&x_{s,\theta_{s-t+1}}
	\end{vmatrix}\hspace{4cm}\nonumber\end{eqnarray}\begin{eqnarray}
	&=&D_{n,t}(i)d_{s,t}(\Theta)-D_{n,s+1}(i)d_{s,t}(\Theta)\nonumber\\ & &+\sum\limits_{j=t}^{s}x_{j,i}\begin{vmatrix}
		x_{t,\theta_1}&x_{t,\theta_2}&...&x_{t,\theta_{s-t+1}}\\
		...&...&...&...\\
		x_{j-1,\theta_1}&x_{j-1,\theta_2}&...&x_{j-1,\theta_{s-t+1}}\\
		D_{n,s+1}(\theta_1)&D_{n,s+1}(\theta_2)&..&D_{n,s+1}(\theta_{s-t+1})\\
		x_{j+1,\theta_1}&x_{j+1,\theta_2}&...&x_{j+1,\theta_{s-t+1}}\\
		...&...&...&...\\
		x_{s,\theta_1}&x_{s,\theta_2}&...&x_{s,\theta_{s-t+1}}
	\end{vmatrix}\nonumber\hspace{5cm}\end{eqnarray}
\begin{eqnarray}
	&=&\sum\limits_{j=s+1}^n(-1)^{s-j}x_{j,i}
	\begin{vmatrix}
		x_{t,\theta_1}&...&x_{t,\theta_{s-t+1}}&0&0&...&0\\
		...&...&...&...&...&...&...\\
		x_{s,\theta_1}&...&x_{s,\theta_{s-t+1}}&0&0&...&0\\
		x_{s+1,\theta_1}&...&x_{s+1,\theta_{s-t+1}}&x_{s+1,s+1}&x_{s+1,s+2}&...&x_{s+1,n}\\
		...&...&...&...&...&...&...\\
		(x_{j,\theta_1}&...&x_{j,\theta_{s-t+1}}&x_{j,s+1}&x_{j,s+2}&...&x_{j,n})\\
		...&...&...&...&...&...&...\\
		x_{n,\theta_1}&...&x_{n,\theta_{s-t+1}}&x_{n,s+1}&x_{n,s+2}&...&x_{n,n}
	\end{vmatrix}\nonumber\\
	&&+D_{n,t}(i)d_{s,t}(\Theta)+\sum\limits_{j=t}^s(-1)^{s-j}
	\begin{vmatrix}
		x_{t,\theta_1}&...&x_{t,\theta_{s-t+1}}&0&...&0\\
		...&...&...&...&...&...\\
		(x_{j,\theta_1}&...&x_{j,\theta_{s-t+1}}&0&...&0)\\
		...&...&...&...&...&...\\
		x_{s,\theta_1}&...&x_{s,\theta_{s-t+1}}&0&...&0\\
		x_{s+1,\theta_1}&...&x_{s+1,\theta_{s-t+1}}&x_{s+1,s+2}&...&x_{s+1,n}\\
		...&...&...&...&...&...\\
		x_{n,\theta_1}&...&x_{n,\theta_{s-t+1}}&x_{n,s+2}&...&x_{n,n}
	\end{vmatrix}\nonumber
	\\
	&=&D_{n,t}(i)d_{s,t}(\Theta)+(-1)^{s-t}
	\begin{vmatrix}
		x_{t,i}&x_{t,\theta_1}&...&x_{t,\theta_{s-t+1}}&x_{t,s+2}&...&x_{t,n}\\
		...&...&...&...&...&...&...\\
		x_{s,i}&x_{s,\theta_1}&...&x_{s,\theta_{s-t+1}}&x_{s,s+2}&...&x_{s,n}\\
		x_{s+1,i}&x_{s+1,\theta_1}&...&x_{s+1,\theta_{s-t+1}}&x_{s+1,s+2}&...&x_{s+1,n}\\
		...&...&...&...&...&...&...\\
		x_{n,i}&x_{n,\theta_1}&...&x_{n,\theta_{s-t+1}}&x_{n,s+2}&...&x_{n,n}
	\end{vmatrix}\nonumber\\ & &\nonumber\\
	&=&D_{n,t}(i)d_{s,t}(\Theta)+(-1)^{s-t}d_{n,t}(i,\theta_1,...,\theta_{s-t+1},s+2,...,n).
	\label{3.30}\end{eqnarray}
Expressions (\ref{3.28})-(\ref{3.30}) imply
\begin{equation}
	\Phi_i(d_{s,t}(\Theta))=D_{n,t}(i)d_{s,t}(\Theta)-d_{n,t}(\tilde{\Theta}_s) \label{3.31}
\end{equation}
where
\begin{equation}
	\tilde{\Theta}_s=\left\lbrace
	\begin{aligned}
		&(0,...,0) &\mathrm{if}~s=n,\\
		&(\theta_1,\theta_2,...,\theta_{s-t+1},i) &\mathrm{if}~s=n-1,\\
		&(\theta_1,\theta_2,...,\theta_{s-t+1},i,s+2,...,n) &\mathrm{if}~s\leqslant n-2.
	\end{aligned}
	\right. \label{3.32}
\end{equation}

For $1\leqslant i_1< i_2<...< i_{n-m+1}\leqslant n$, $1\leqslant m\leqslant n-1$ and $f\in\msr C_n$, we define
\begin{equation}
	\begin{bmatrix}
		E_{i_1,n+1}&E_{i_2,n+1}&...&E_{i_{n-m+1},n+1}\\
		x_{m,i_1}&x_{m,i_2}&...&x_{m,i_{m-n+1}}\\
		x_{m+1,i_1}&x_{m+1,i_2}&...&x_{m+1,i_{m-n+1}}\\
		...&...&...&...\\
		x_{n-1,i_1}&x_{n-1,i_2}&...&x_{n-1,i_{m-n+1}}
	\end{bmatrix}(f)=
	\sum\limits_{j=1}^{n-m+1}E_{i_j,n+1}(a_{i_j}f), \label{3.33}
\end{equation}
where
\begin{equation}a_{i_j}=(-1)^{1+j}d_{n-1,m}(i_1,...,i_{j-1},i_{j+1},...,i_{n-m+1})\end{equation}is the algebraic cofactor of $E_{i_j,n+1}$ in the matrix (cf. (2.48)).
By (2.84), (\ref{3.33}) is equal to
\begin{equation}
	\sum\limits_{j=1}^{n-m+1}\left(\Phi_{i_j}(a_{i_j}f)-\Psi_{i_j}(a_{i_j}f) \right)
	=\sum\limits_{j=1}^{n-m+1}\left(\Phi_{i_j}(a_{i_j})f+a_{i_j}\Phi_{i_j}(f)-\Psi_{i_j}(a_{i_j}f) \right) \label{3.34}
\end{equation}

According to (\ref{3.31}),
\begin{eqnarray}
	\sum\limits_{j=1}^{n-m+1}\Phi_{i_j}(a_{i_j})&=&
	\sum\limits_{j=1}^{n-m+1}\left( D_{n,m}(i_j)a_{i_j}-(-1)^{n-m}d_{n,m}(i_1,...,i_{n-m+1})\right) \nonumber\\
	&=&\sum\limits_{j=1}^{n-m+1}\left( x_{n,i_{j}}a_{i_j}+\sum\limits_{t=m}^{n-1}x_{t,i_j}D_{n,t+1}(t+1)a_{i_j}\right) \nonumber\\& &-(-1)^{n-m}(n-m+1)d_{n,m}(i_1,...,i_{n-m+1})\nonumber\\
	&=&\sum\limits_{j=1}^{n-m+1}x_{n,i_j}a_{i_j}+\sum\limits_{t=m}^{n-1}\left( D_{n,t+1}(t+1)\sum\limits_{j=1}^{n-m+1}x_{t,i_j}a_{i_j}\right)\nonumber\\ & &-(-1)^{n-m}(n-m+1)d_{n,m}(i_1,...,i_{n-m+1})\nonumber\\
	&=&\!\!(-1)^{n-m}d_{n,m}(i_1,...,i_{n-m+1})-(-1)^{n-m}(n-m+1)d_{n,m}(i_1,...,i_{n-m+1})\nonumber\\
	&=&(m-n)(-1)^{n-m}d_{n,m}(i_1,...,i_{n-m+1}), \label{3.35}
\end{eqnarray}
where the second equality was obtained by Lemma 2.2 (ii). Moreover,
\begin{eqnarray}
	\sum\limits_{j=1}^{n-m+1}a_{i_j}\Phi_{i_j}&=&\sum\limits_{1\leqslant t\leqslant s\leqslant n}\sum\limits_{j=1}^{n-m+1}a_{i_j}D_{n,s}(t)x_{s,i_j}\partial_{x_{s,t}}
	\nonumber\\ &=&\sum\limits_{1\leqslant t\leqslant s\leqslant n}D_{n,s}(t)
	\begin{vmatrix}
		x_{s,i_1}&x_{s,i_2}&...&x_{s,i_{n-m+1}}\\
		x_{m,i_1}&x_{m,i_2}&...&x_{m,i_{n-m+1}}\\
		x_{m+1,i_1}&x_{m+1,i_2}&...&x_{m+1,i_{n-m+1}}\\
		...&...&...&...\\
		x_{n-1,i_1}&x_{n-1,i_2}&...&x_{n-1,i_{n-m+1}}
	\end{vmatrix}\partial_{x_{s,t}}\nonumber\\
	&=&(-1)^{n-m}d_{n,m}(i_1,...,i_{n-m+1})\sum\limits_{t=1}^{n}x_{n,t}\partial_{x_{n,t}}\nonumber\\ & &+\sum\limits_{1\leqslant t\leqslant s<m}D_{n,s}(t)
	\begin{vmatrix}
		x_{s,i_1}&x_{s,i_2}&...&x_{s,i_{n-m+1}}\\
		x_{m,i_1}&x_{m,i_2}&...&x_{m,i_{n-m+1}}\\
		x_{m+1,i_1}&x_{m+1,i_2}&...&x_{m+1,i_{n-m+1}}\\
		...&...&...&...\\
		x_{n-1,i_1}&x_{n-1,i_2}&...&x_{n-1,i_{n-m+1}}
	\end{vmatrix}\partial_{x_{s,t}}.
	\label{3.36}
\end{eqnarray}
Furthermore,
\begin{eqnarray}	\!\!\!\!\!\!\!\!\!\!\!\!\!\!\!\sum\limits_{j=1}^{n-m+1}\Psi_{i_j}a_{i_j}\!\!\!\!\!&=&\!\!\!\!\!\sum\limits_{j=1}^{n-m+1}\sum\limits_{s=1}^nk_{n-s+1}D_{n,s}(i_j)a_{i_j}\nonumber\\		\!\!\!\!\!&=&\!\!\!\!\!\sum\limits_{s=m}^n\sum\limits_{j=1}^{n-m+1}k_{n-s+1}D_{n,s}(i_j)a_{i_j}
	+\sum\limits_{s=1}^{m-1}\sum\limits_{j=1}^{n-m+1}k_{n-s+1}D_{n,s}(i_j)a_{i_j}\nonumber\\		\!\!\!\!\!\!\!\!\!\!&=&\!\!\!\!\!\sum\limits_{s=m}^n(-1)^{n-m}k_{n-s+1}d_{n,m}(i_1,...,i_{n-m+1})
	+\sum\limits_{s=1}^{m-1}\sum\limits_{j=1}^{n-m+1}k_{n-s+1}D_{n,s}(i_j)a_{i_j}.
	\label{3.37}
\end{eqnarray}
In summary, we have
\begin{eqnarray}
	&&\begin{bmatrix}
		E_{i_1,n+1}&E_{i_2,n+1}&...&E_{i_{n-m+1},n+1}\\
		x_{m,i_1}&x_{m,i_2}&...&x_{m,i_{m-n+1}}\\
		x_{m+1,i_1}&x_{m+1,i_2}&...&x_{m+1,i_{m-n+1}}\\
		...&...&...&...\\
		x_{n-1,i_1}&x_{n-1,i_2}&...&x_{n-1,i_{m-n+1}}
	\end{bmatrix}\nonumber\\
	&=&(-1)^{n-m}d_{n,m}(i_1,...,i_{n-m+1})\left(m-n+\sum\limits_{t=1}^nx_{n,t}\partial_{x_{n,t}}-\sum\limits_{s=1}^{n-m+1}k_s \right)\nonumber\\
	&&+\sum\limits_{1\leqslant t\leqslant s<m}\sum\limits_{j=1}^{n-m+1}a_{i_j}D_{n,s}(t)x_{s,i_j}\partial_{x_{s,t}}
-\sum\limits_{s=1}^{m-1}\sum\limits_{j=1}^{n-m+1}k_{n-s+1}D_{n,s}(i_j)a_{i_j}.
	\label{3.38}
\end{eqnarray}

For $1\leqslant j\leqslant i\leqslant n$, we set
\begin{equation}\mbb J_{i,j}=\{\Theta=(\sta_1,\sta_2,...,\sta_{i-j+1})\mid \sta_r\in\ol{1,n};\;d_{i,j}(\Theta)\neq 0\}.\end{equation}
Recall that $\lmd=\sum_{i=1}^nk_i\lmd_i$ is a dominant integral weight of $sl(n+1)$. Denote
\begin{eqnarray}
	S_n(\lambda)&=&\big\{\prod\limits_{1\leqslant j\leqslant i\leqslant n}\prod\limits_{\Theta\in \mbb J_{i,j}} d_{i,j}(\Theta)^{\alpha_{i,j}(\Theta)} \mid \alpha_{i,j}(\Theta)\in\mbb{N};\nonumber\\ & &
	\sum\limits_{i=j}^n\sum\limits_{\Theta\in\mbb J_{i,j}}\alpha_{i,j}(\Theta)\leqslant k_{n-j+1}~\mathrm{for}~j\in\overline{1,n}\big\}\subset\msr C_n.  \label{3.39}
\end{eqnarray}
Moreover,
the homogeneous subset of $S_n(\lambda)$ with degree r in $\{x_{n,1},x_{n,2},...x_{n,n}\}$ is
\begin{equation}
	(S_n(\lambda))_r=\big\lbrace \prod\limits_{1\leqslant j\leqslant i\leqslant n}\prod\limits_{\Theta\in\mbb J_{i,j}}d_{i,j}(\Theta)^{\alpha_{i,j}(\Theta)}\in S_n(\lambda)~\Bigg|~\sum\limits_{j=1}^n\sum\limits_{\Theta\in\mbb J_{n,j}}\alpha_{n,j}(\Theta)=r\big\rbrace  \label{3.40}
\end{equation}
Recall $|\lmd|=\sum_{i=1}^nk_i$. Then
\begin{equation}
	S_n(\lambda)=\bigcup\limits_{r=0}^{|\lmd|}(S_n(\lambda))_r \label{3.41}
\end{equation}
According to (1.10)-(1.12) and (2.56), the irreducible $sl(n+1)$-module
\begin{equation}V_n(\lmd)=\ol{M}\subset\wht M\subset\msr C_n.\end{equation}

{\bf Theorem 2.4}\quad {\it
	As a vector space, $V_n(\lambda)$ is spanned by $S_n(\lambda)$. Moreover, the homogeneous subspace $(V_n(\lambda))_r=\ol M_r$ is spanned by $(S_n(\lambda))_r$.}

\begin{proof}
	We prove the proposition by induction on $n$. If $n=2$, the sets
	\begin{equation}\mbb J_{1,1}=\mbb J_{2,2}=\{(1),(2)\},\qquad \mbb J_{2,1}=\{(1,2),(2,1)\}.\end{equation}
	Moreover,
	\begin{equation}d_{1,1}(1)=x_{1,1},\;\;d_{1,1}(2)=-1,\;\;d_{2,2}(1)=x_{2,1},\;\;d_{2,2}(2)=x_{2,2},\end{equation}
	\begin{equation}d_{2,1}(1,2)=-d_{2,1}(2,1)=\left|\begin{array}{cc}x_{1,1}&-1\\ x_{2,1}&
			x_{2,2}\end{array}\right|=x_{1,1}x_{2,2}+x_{2,1}.\end{equation}
	Under the identification (2.53), the basis (2.27) is contained in the set
	\begin{eqnarray}S_2(\lambda)&\subset&
		\big\{ \pm x_{2,1}^{\alpha_{2,1}}x_{2,2}^{\alpha_{2,2}}x_{1,1}^{\alpha_{1,1}}(x_{2,1}+x_{2,2}x_{1,1})^{\alpha'_{1,1}} \nonumber\\ & &\mid 0\leqslant\alpha_{2,1}+\alpha_{2,2}\leqslant k_1~\mathrm{and}~0\leqslant\alpha_{1,1}+\alpha'_{1,1}\leqslant k_2 \big\} \subset V(\lmd)\label{3,42}
	\end{eqnarray}
	(cf. (2.18) and (2.19)). Thus $V_2(\lmd)=\ol M=\mbox{Span}\:S_2(\lmd)$ and
	$(V_2(\lmd))_r=\ol M_r=\mbox{Span}\:(S_2(\lmd))_r$. So the theorem holds for $n=2$.\psp
	
	Set
	\begin{equation}W(\lambda)_r=\mbox{Span}\:(S_n(\lmd))_r.\end{equation}
	Suppose that the theorem  holds for $n-1$, which is equivalent to
	\begin{equation} (V_n(\lambda))_0=V_{n-1}(k_2\lmd_1+k_3\lmd_2+\cdots+k_n\lmd_{n-1})=W(\lambda)_0.\end{equation}
	Assume $W(\lambda)_r=(V_n(\lambda))_r$ for some $0\leqslant r<|\lmd|=\sum\limits_{i=1}^{n}k_i$. Take $f=\prod\limits_{1\leqslant j\leqslant i\leqslant n}\prod\limits_{\Theta\in \mbb J_{i,j}}d_{i,j}(\Theta)^{\alpha_{i,j}(\Theta)}\in (S_n(\lambda))_r$. By (\ref{3.31}) and
	(\ref{3.32}),
	\begin{eqnarray}
		E_{s,n+1}(f)&=&(\Phi_s-\Psi_s)\left(\prod\limits_{1\leqslant j\leqslant i\leqslant n}\prod\limits_{\Theta\in \mbb J_{i,j}}d_{i,j}(\Theta)^{\alpha_{i,j}(\Theta)} \right) \nonumber\\
		&=&f \sum\limits_{1\leqslant j\leqslant i\leqslant n}\sum\limits_{\Theta\in\mbb J_{i,j}}\alpha_{i,j}(\Theta)d_{i,j}(\Theta)^{-1}\left( D_{n,j}(s)d_{i,j}(\Theta)-d_{n,j}(\tilde{\Theta}_s)\right)\nonumber\\ & &-f\sum\limits_{j=s}^nk_{n-j+1}D_{n,j}(s)
		\nonumber\\ &=&\sum\limits_{j=s}^n\left( \sum\limits_{i=j}^n\sum\limits_{\Theta\in\mbb J_{i,j}}\alpha_{i,j}(\Theta)-k_{n-j+1}\right)D_{n,j}(s)f\nonumber\\& &-\sum\limits_{1\leqslant j\leqslant i\leqslant n}\sum\limits_{\Theta\in\mbb J_{i,j}}\alpha_{i,j}(\Theta)d_{i,j}(\Theta)^{-1}d_{n,j}(\tilde{\Theta}_s)f\in W(\lmd)_{r+1},		 \label{3.45}
	\end{eqnarray}
	where we have used the fact $D_{n,j}(s)=0$ when $j<s$, and $D_{n,j}(s)f\in (S_n(\lambda))_{r+1}$ if $ \sum\limits_{i=j}^n\sum\limits_{\Theta\in\mbb J_{i,j}}\alpha_{i,j}(\Theta)<k_{n-j+1}$. Therefore
	\begin{equation}(V_n(\lambda))_{r+1}=\sum\limits_{s=1}^nE_{s,n+1}((V_n(\lambda))_r)\subseteq W(\lambda)_{r+1}
	\end{equation}by (2.9) and (2.13).
	
	On the other hand, we set
	\begin{equation}
		P_{r+1}^m=(S_n(\lambda))_{r+1}\bigcap\big\lbrace \prod\limits_{1\leqslant j\leqslant i\leqslant n}\prod\limits_{\Theta\in\mbb J_{i,j}}d_{i,j}(\Theta)^{\alpha_{i,j}(\Theta)}~\Bigg|~\sum\limits_{j=1}^{m}\sum\limits_{\Theta\in\mbb J_{n,j}}\alpha_{n,j}(\Theta)> 0 \big\rbrace
	\end{equation}
	for $1\leqslant m\leqslant n$. It is obvious that
	\begin{equation}
		P_{r+1}^1\subseteq P_{r+1}^2\subseteq\cdots\subseteq P_{r+1}^n=(S_n(\lambda))_{r+1}.\end{equation}We may assume $\sum\limits_{s=1}^{n'-1}k_s\leqslant r<\sum\limits_{s=1}^{n'}k_s$ with some $n'\in\overline{1,n}$.
	We claim
	\begin{equation}
		P_{r+1}^{n-n'+1}=(S_n(\lambda))_{r+1}.\label{2.112}
	\end{equation}
	
	Assume that it is not true. Take any
	\begin{equation}   \prod\limits_{1\leqslant j\leqslant i\leqslant n}\prod\limits_{\Theta\in\mathbb{J}_{i,j}}d_{i,j}(\Theta)^{\alpha_{i,j}(\Theta)}\in (S_n(\lambda))_{r+1}\setminus P_{r+1}^{n-n'+1}. \end{equation}
	Then
	\begin{equation}
		\sum\limits_{j=1}^{n-n'+1}\sum\limits_{\Theta\in\mathbb{J}_{n,j}}\alpha_{n,j}(\Theta)=0,
	\end{equation}
	and
	\begin{equation}
		\sum\limits_{j=n-n'+2}^{n}\sum\limits_{\Theta\in\mathbb{J}_{n,j}}\alpha_{n,j}(\Theta)=r+1>\sum\limits_{s=1}^{n'-1}k_s.
	\end{equation}
	However, by (2.99),
	\begin{equation}
		\sum\limits_{j=n-n'+2}^{n}\sum\limits_{\Theta\in\mathbb{J}_{n,j}}\alpha_{n,j}(\Theta)\leqslant
		\sum\limits_{j=n-n'+2}^n\sum\limits_{i=j}^n\sum\limits_{\Theta\in\mathbb{J}_{i,j}}\alpha_{i,j}(\Theta)\leqslant
		\sum\limits_{j=n-n'+2}^nk_{n-j+1}=\sum\limits_{s=1}^{n'-1}k_s.
	\end{equation}
	which contradicts  (2.116). Thus
	\begin{equation}
		P_{r+1}^1\subseteq P_{r+1}^2\subseteq\cdots\subseteq P_{r+1}^{n-n'+1}=(S_n(\lambda))_{r+1}.
	\end{equation}
	
	Take any $g\in P_{r+1}^1$. It can be written as $g=g'D_{n,1}(1)$ for some $g'\in (S_n(\lambda))_r$. By the definition of $S_n(\lambda)$,
	\begin{equation}
		g'd_{n-1,1}(\Theta)\in (S_n(\lambda))_r\subseteq (V_n(\lambda))_r
	\end{equation}
	for any $\Theta\in\mbb J_{n-1,1}$. According to (2.85), (2.86), (2.91) and (2.107),
	\begin{equation}
		\begin{bmatrix}
			E_{1,n+1}&E_{2,n+1}&...&E_{n,n+1}\\
			x_{1,1}&x_{1,2}&...&x_{1,n}\\
			x_{2,1}&x_{2,2}&...&x_{2,n}\\
			...&...&...&...\\
			x_{n-1,1}&x_{n-1,2}&...&x_{n-1,n}
		\end{bmatrix}	\left( g' \right) =(-1)^{n-1}\left(r-\sum\limits_{i=1}^nk_i-n+1 \right)g.
	\end{equation}
	Since the coefficient $r-\sum\limits_{i=1}^nk_i-n+1< -n+1$ is nonzero,  $g\in (V_n(\lambda))_{r+1}$. Hence $P_{r+1}^1\subseteq (V_n(\lambda))_{r+1}$.\pse

	Next we assume that $P_{r+1}^{m-1}\subseteq (V_n(\lambda))_{r+1}$ for some $m\in\overline{2,n-n'+1}$ by (2.118). Take any  $h\in P_{r+1}^m\backslash P_{r+1}^{m-1}$. Without loss of generality, we assume
	$h=d_{n,m}(\vartheta_1,...,\vartheta_{n-m+1})h'$
	with $h'\in (S_n(\lambda))_r$ and $1\leqslant \vartheta_1<\cdots<\vartheta_{n-m+1}\leqslant n$.
	According to (2.98),
	\begin{equation}
		h'd_{n-1,m}(\Theta)\in (S_n(\lambda))_r\qquad\mbox{for any}\;\; \Theta\in\mbb J_{n-1,m}. \end{equation}
	Suppose
	\begin{equation} h'=\prod\limits_{1\leqslant j\leqslant i\leqslant n}\prod\limits_{\Theta\in\mbb J_{i,j}}d_{i,j}(\Theta)^{\beta_{i,j}(\Theta)}.\end{equation}
	We write $h'=h'_1h'_2h'_3$ according to the range of the index $j$ as follows:
	\begin{equation}h'_1=\prod\limits_{j=1}^{m-1}\prod\limits_{i=j}^{n-1}\prod\limits_{\Theta\in \mbb J_{i,j}}d_{i,j}(\Theta)^{\beta_{i,j}(\Theta)},\;\;h'_2=\prod\limits_{j=m}^{n-1}\prod\limits_{i=j}^{n-1}\prod\limits_{\Theta\in\mbb J_{i,j}}d_{i,j}(\Theta)^{\beta_{i,j}(\Theta)},\end{equation}
	\begin{equation}
		h'_3=\prod\limits_{j=m}^{n}\prod\limits_{\Theta\in\mbb J_{n,j}}d_{n,j}(\Theta)^{\beta_{n,j}(\Theta)}
	\end{equation}
	based on (2.111) and $h\not\in P_{r+1}^{m-1}$.
	We apply the operator defined in (\ref{3.33}). By (\ref{3.34})-(\ref{3.38}),
	\begin{eqnarray}
		& &\begin{bmatrix}
			E_{\vartheta_1,n+1}&E_{\vartheta_2,n+1}&...&E_{_{n-m+1},n+1}\\
			x_{m,\vartheta_1}&x_{m,\vartheta_2}&...&x_{m,\vartheta_{m-n+1}}\\
			x_{m+1,\vartheta_1}&x_{m+1,\vartheta_2}&...&x_{m+1,\vartheta_{m-n+1}}\\
			...&...&...&...\\
			x_{n-1,\vartheta_1}&x_{n-1,\vartheta_2}&...&x_{n-1,\vartheta_{m-n+1}}
		\end{bmatrix}(h')\nonumber
		\\	&=&(-1)^{n-m}(m-n-\sum\limits_{s=1}^{n-m+1}k_s)h'd_{n,m}(\vartheta_1,...,\vartheta_{n-m+1}\nonumber)\\ &&+\sum\limits_{l=1}^{n-m+1}a_{\vartheta_{l}}\Phi_{\vartheta_l}(h')
		-\sum\limits_{l=1}^{n-m+1}\sum\limits_{s=1}^{m-1}k_{n-s+1}D_{n,s}(\vartheta_l)a_{\vartheta_l}h'\nonumber\\
		&=&(-1)^{n-m}(m-n-\sum\limits_{s=1}^{n-m+1}k_s)h-\sum\limits_{l=1}^{n-m+1}\sum\limits_{s=1}^{m-1}k_{n-s+1}D_{n,s}
		(\vartheta_l)a_{\vartheta_l}h'\nonumber\\
		&&+\sum\limits_{l=1}^{n-m+1}a_{\vartheta_{l}}\left( \Phi_{\vartheta_l}(h'_1)h'_2h'_3+\Phi_{\vartheta_l}(h'_2)h'_1h'_3+\Phi_{\vartheta_l}(h'_3)h'_1h'_2\right).
		\label{3.54}
	\end{eqnarray}
	The first term is an integral multiple of $h$. Let us figure out what the remnant is. As (\ref{3.45}), we have
	\begin{eqnarray}& & \sum\limits_{l=1}^{n-m+1}a_{\vartheta_l}\Phi_{\vartheta_l}(h'_1)h'_2h'_3
		-\sum\limits_{l=1}^{n-m+1}\sum\limits_{s=1}^{m-1}k_{n-s+1}D_{n,s}(\vartheta_l)a_{\vartheta_l}h'\nonumber\\
		&=&\sum\limits_{l=1}^{n-m+1}a_{\vartheta_l}h'_2h'_3\left( \Phi_{\vartheta_l}(h'_1)-\sum\limits_{s=1}^{m-1}k_{n-s+1}D_{n,s}(\vartheta_l)h'_1   \right)  \label{3.55}
		\nonumber\\  		&=&\sum\limits_{l=1}^{n-m+1}a_{\vartheta_l}h'_2h'_3\big[  \sum\limits_{j=1}^{m-1}\sum\limits_{i=1}^n\sum\limits_{\Theta\in\mathbb{J}_{i,j}}\beta_{i,j
		}(\Theta)d_{i,j}(\Theta)^{-1}\big( D_{n,j}(\vartheta_l)d_{i,j}(\Theta)-d_{n,j}(\tilde{\Theta}_{\vartheta_l})
		\big)h'_1  \nonumber\\& &-\sum\limits_{s=1}^{m-1}k_{n-s+1}D_{n,s}(\vartheta_{l})h'_1\big]\nonumber\\
		&=&\sum\limits_{l=1}^{n-m+1}\sum\limits_{j=1}^{m-1}\big[ \big( \sum\limits_{i=1}^{n}\sum\limits_{\Theta\in\mathbb{J}_{i,j}}\beta_{i,j}(\Theta)-k_{n-j+1}\big) D_{n,j}(\vartheta_l)h'a_{\vartheta_{l}} \nonumber\\& &-\sum\limits_{i=1}^n\sum\limits_{\Theta\in\mathbb{J}_{i,j}}\beta_{i,j}(\Theta)d_{i,j}(\Theta)^{-1}d_{n,j}
		(\tilde{\Theta}_{\vartheta_l})h'a_{\vartheta_l}\big].
		\label{3.56}
	\end{eqnarray}
	Since $j\leqslant m-1$,  $\beta_{i,j}(\Theta)\neq 0$ implies \begin{equation}d_{i,j}(\Theta)^{-1}d_{n,j}(\tilde{\Theta}_{\vartheta_l})h'a_{\vartheta_{l}}\in P_{r+1}^{m-1}. \end{equation}
	If $\sum\limits_{i=1}^{n}\sum\limits_{\Theta\in\mbb J_{i,j}}\beta_{i,j}(\Theta)-k_{n-j+1}\neq 0$, (\ref{3.39}) yields $D_{n,j}(\vartheta_l)h'a_{\vartheta_l}\in P_{r+1}^{m-1}$. Thus  ($\ref{3.55}$) is a polynomial in $\mbox{Span}\:P_{r+1}^{m-1}$.
	
	Note that $h'_2$ is a polynomial in $\{x_{i,j}\mid m\leqslant j\leqslant i\leqslant n-1\}$. By (\ref{3.36})
	\begin{eqnarray} 	    		&&\sum\limits_{l=1}^{n-m+1}a_{\vartheta_l}\Phi_{\vartheta_l}(h'_2)h'_1h'_3\nonumber\\ &=&(-1)^{n-m}h'_1h'_3d_{n,m}(\vartheta_1,...,\vartheta_{n-m+1})\sum\limits_{t=1}^{n}x_{n,t}\partial_{x_{n,t}}(h'_2)\nonumber\\
		&&+h'_1h'_3\!\!\!\sum\limits_{1\leqslant t\leqslant s<m}\!\!\!D_{n,s}(t)
		\begin{vmatrix}
			x_{s,\vartheta_1}&x_{s,\vartheta_2}&...&x_{s,\vartheta_{n-m+1}}\\
			x_{m,\vartheta_1}&x_{m,\vartheta_2}&...&x_{m,\vartheta_{n-m+1}}\\
			x_{m+1,\vartheta_1}&x_{m+1,\vartheta_2}&...&x_{m+1,\vartheta_{n-m+1}}\\
			...&...&...&...\\
			x_{n-1,\vartheta_1}&x_{n-1,\vartheta_2}&...&x_{n-1,\vartheta_{n-m+1}}
		\end{vmatrix}
		\partial_{x_{s,t}}(h'_2)=0
		\label{3.57}
	\end{eqnarray}
	Now there leaves one term in (\ref{3.54}). Since $h'_3$ is a polynomial of $\{x_{i,j}\mid m\leqslant j\leqslant i\leqslant n\}$ and is homogeneous in $X_n$ with degree $r$, we have
	\begin{eqnarray}    \sum\limits_{l=1}^{n-m+1}a_{\vartheta_l}\Phi_{\vartheta_l}(h'_3)h'_1h'_2 &=&\!(-1)^{n-m}h'_1h'_2d_{n,m}(\vartheta_1,...,\vartheta_{n-m+1})\sum\limits_{t=1}^{n}x_{n,t}\partial_{x_{n,t}}(h'_3)\nonumber\\
		&+&\!\!\!\!h'_1h'_2\!\!\sum\limits_{1\leqslant t\leqslant s<m}\!\!D_{n,s}(t)
		\begin{vmatrix}
			x_{s,\vartheta_1}&x_{s,\vartheta_2}&...&x_{s,\vartheta_{n-m+1}}\\
			x_{m,\vartheta_1}&x_{m,\vartheta_2}&...&x_{m,\vartheta_{n-m+1}}\\
			x_{m+1,\vartheta_1}&x_{m+1,\vartheta_2}&...&x_{m+1,\vartheta_{n-m+1}}\\
			...&...&...&...\\
			x_{n-1,\vartheta_1}&x_{n-1,\vartheta_2}&...&x_{n-1,\vartheta_{n-m+1}}
		\end{vmatrix}
		\partial_{x_{s,t}}(h'_3)\nonumber\\
		&=&\!(-1)^{n-m}rh'_1h'_2h'_3d_{n,m}(\vartheta_1,...,\vartheta_{n-m+1})\nonumber\\
		&=&\!(-1)^{n-m}rh.
		\label{3.58}
	\end{eqnarray}
	Therefore, by (\ref{3.55})-(\ref{3.58}),
	\begin{equation}
		(\ref{3.54})=(-1)^{n-m}(r+m-n-\sum\limits_{s=1}^{n-m+1}k_s)h+h'',
	\end{equation}
	with some $h''\in \mbox{Span}\: P_{r+1}^m$. Since $m\leqslant n-n'+1$, we have $r+m-n\leqslant r<\sum\limits_{s=1}^{n'}k_s\leqslant\sum\limits_{s=1}^{n-m+1}k_s$ and so the coefficient $(-1)^{n-m}(r+m-n-\sum\limits_{s=1}^{n-m+1}k_s)$ is nonzero. By ($\ref{2.112}$), $h$ lies in $V(\lambda)_{r+1}$. Based on  induction on $m$, we have $P_{r+1}^m\subseteq (V_n(\lambda))_{r+1}$ for all $1\leqslant m\leqslant n-n'+1$. Hence $(S_n(\lambda))_{r+1}\subseteq (V_n(\lambda))_{r+1}$ and $(V_n(\lambda))_{r+1}=W(\lambda)_{r+1}$ (cf. (2.107)).

	It leaves to show $(V_n(\lambda))_r=\left\lbrace 0\right\rbrace $ when $r>\sum\limits_{i=1}^{n}k_i$. Recall $|\lmd|=\sum_{i=1}^nk_i$. If $f'=\prod\limits_{1\leqslant j\leqslant i\leqslant n}\prod\limits_{\Theta\in\mbb J_{i,j}}d_{i,j}(\Theta)^{\gamma_{i,j}(\Theta)}$ is any nonzero polynomial in $(S_n(\lambda))_{|\lmd|}$. Since $\sum\limits_{j=1}^n\sum\limits_{\Theta\in\mbb J_{n,j}}\gamma_{n,j}(\Theta)=\sum\limits_{i=1}^nk_i$,
	(\ref{3.39}) shows
	\begin{equation}
		f'=\prod\limits_{j=1}^n\prod\limits_{\Theta\in\mbb J_{n,j}}d_{n,j}(\Theta)^{\gamma_{n,j}(\Theta)}
	\end{equation}
	with $\sum\limits_{\Theta\in\mbb J_{n,j}}\gamma_{n,j}(\Theta)=k_{n-j+1}$ for $1\leqslant j\leqslant n$. As (\ref{3.45}), we have
	\begin{equation}
		E_{s,n+1}(f)=\sum\limits_{j=s}^n\left( \sum\limits_{\Theta\in\mbb J_{n,j}}\gamma_{n,j}(\Theta)-k_{n-j+1} \right) D_{n,j}(s)f'=0
	\end{equation}
	for $1\leqslant s\leqslant n$. Hence $(V_n(\lambda))_{|\lmd|+1}=\{0\}$. Thereby, $(V_n(\lambda))_{r}=\{0\}$ for any
	$r>|\lmd|$ by (2.13).
\end{proof}

\section{Bases and Singular Vectors}

 In this section, we determine a basis for the $sl(n+1)$-module $V_n(\lmd)$ realized in Theorem 2.4 and a connection with Gelfand-Tsetlin basis (cf. \cite{GT1}). Moreover, we calculate all the $sl(n)$-singular vectors in $\wht M$. In particular, the results  restricted to $\wht M_r$ and $V_n(\lmd)=\ol M$ yield to certain sum-product type combinatorial identities.

 First we define a lexicographic degree \begin{equation}\mathrm{lexdeg}(X^{\alpha})=(\alpha_{1,1},\alpha_{2,1},\alpha_{2,2},...,\alpha_{j,1},\alpha_{j,2}...,
 \alpha_{j,j},...,\alpha_{n,n})\end{equation}
 for monomials $\mathscr{C}_n$ according to the order:
\begin{equation}
	(1,1)<(2,1)<(2,2)\cdots<(j,1)<(j,2)\cdots <(j,j)\cdots <(n,n).
\end{equation}
We  say that $\mathrm{lexdeg}(X^{\alpha})>\mathrm{lexdeg}(X^{\alpha'})$ if there exists $(s,t)$ in (3.2)  such that $\alpha_{i,j}=\alpha'_{i,j}$ for all indexes $(i,j)<(s,t)$  and $\alpha_{s,t}>\alpha'_{s,t}$. This gives a total order. For any polynomial $f\in\mathscr{C}_n$, we define $\mathrm{lexdeg}(f)$ to be the maximal lexicographic degree of monomial appears in $f$. Define a relation in $(V_n(\lambda))_r$ ($0\leqslant r\leqslant|\lmd|$) as follows:
\begin{equation}
	\forall u,v\in (V_n(\lambda))_r,~u\sim v~~~\mathrm{if}~~\mathrm{lexdeg}(u)=\mathrm{lexdeg}(v).
\end{equation}
 This is an equivalence relation. Let
 \begin{equation}{\cal B}(\lambda)_r\;\;\mbox{ be a set of representatives of the elements in}\;\; (S_n(\lambda))_r/\sim\end{equation}
  (cf. (2.99) and (2.100)), and $\mathcal{Q}(\lambda)_r$ be the set of monomials with maximal lexicographic degree of the elements in $\mathcal{B}(\lambda)_r$.
 Suppose
 \begin{equation}f=\prod\limits_{1\leqslant j\leqslant i\leqslant n}\prod\limits_{\Theta\in\mathbb J_{i,j}}d_{i,j}(\Theta)^{\alpha_{i,j}(\Theta)}\in (S_n(\lambda))_r. \end{equation}
  Without loss of generality, we can rearrange the index $\Theta=(\theta_1,...,\theta_{i-j+1})\in\mathbb J_{i,j}$ (cf. (2.98)) in each factor $d_{i,j}(\Theta)$ of $f$ such that $\theta_1<...<\theta_{i-j+1}$. Thus the monomial with maximal lexcoghraphic degree in $d_{i,j}(\Theta)$ is $\prod\limits_{s=1}^{i-j+1}x_{j+s-1,\theta_s}$. Therefore, the the monomial with maximal lexcoghraphic degree in $f$ is
\begin{equation}
	\prod\limits_{1\leqslant j\leqslant i\leqslant n}\prod\limits_{\substack{\Theta\in \mathbb J_{i,j}\\ \theta_1<\cdots<\theta_{i-j+1}}}\prod\limits_{s=1}^{i-j+1}x_{j+s-1,\theta_s}^{\alpha_{i,j}(\Theta)}.
\end{equation}

By the definition of $S_n(\lambda)$ in (2.99) and (3.6), we conclude that $\mathcal{Q}(\lambda)_r$ is the set of monomials $X^{\alpha}\in \mathscr{C}_n$ satisfying the following conditions:
\begin{equation}
	\sum\limits_{j=1}^{n}\alpha_{n,j}=r
\end{equation}
and
\begin{equation}
	\begin{aligned}
		0\leqslant&\alpha_{n,1}\leqslant k_1, \\
		0\leqslant&\alpha_{n,1}+\alpha_{n,2}\leqslant \alpha_{n-1,1}+k_1\leqslant k_1+k_2, \\
		&...\\
		0\leqslant&\sum\limits_{j=1}^{l}\alpha_{n,j}\leqslant \sum\limits_{j=1}^{l-1}\alpha_{n-1,j}+k_1\leqslant\cdots\leqslant\sum\limits_{j=1}^{l-s}\alpha_{n-s,j}
+\sum\limits_{i=1}^{s}k_i\leqslant\cdots\leqslant \sum\limits_{i=1}^{l}k_i, \\
		&...\\
		0\leqslant&\sum\limits_{j=1}^{n}\alpha_{n,j}\leqslant\sum\limits_{j=1}^{n-1}\alpha_{n-1,j}+k_{1}\leqslant\cdots
\leqslant\alpha_{1,1}+\sum\limits_{i=1}^{n-1}k_i\leqslant\sum\limits_{i=1}^{n}k_i.
	\end{aligned}
\end{equation}

By Theorem 2.4,	
\begin{equation}\mathcal{B}(\lambda)=\bigcup\limits_r\mathcal{B}(\lambda)_r\end{equation}
 is the polynomial basis of $V_n(\lambda)$. Set
 \begin{equation}\mathcal{Q}(\lambda)=\bigcup\limits_r\mathcal{Q}(\lambda)_r. \end{equation}
 Next we associate $X^{\alpha}\in\mathcal{Q}(\lambda)$ to a Gelfand-Tsetlin pattern of shape $\lambda$. By Theorem 2.3, let \begin{equation}\tau_{n+1,i}=E_{n+2-i,n+2-i}(1)=\frac{1}{n+1}c_{n+1}-c_n-\sum\limits_{s=1}^{i-1}k_s\;\; \for \;\;i\in\overline{2,n+1}, \end{equation} and set
 \begin{equation}\tau_{n+1,1}=\frac{1}{n+1}c_{n+1}-c_n.\end{equation}
Put
\begin{equation}\ves_{i}=(0,...,0,\stl{i}{1},0,...0).\end{equation}
Then we can extend the weight $\lmd$  of $sl(n+1)$ to that of $gl(n+1)$
\begin{equation}
	\lambda=\tau_{n+1,1}\varepsilon_1+\tau_{n+1,2}\varepsilon_2+\cdots+\tau_{n+1,n+1}\varepsilon_{n+1}
\end{equation}
(cf. (2.1), (2.8), (2.73) and (2.74)). Define
\begin{equation}\tau_{i,j}=\tau_{i+1,j+1}+\alpha_{j-i+n,n-i+1}\;\;\for\;\;1\leqslant j \leqslant i\leqslant n\end{equation}
 inductively.
 Expression (3.8) yields
\begin{equation}
	\tau_{i+1,j}\geqslant \tau_{i,j}\geqslant\tau_{i+1,j+1}\qquad\for\;\; 1\leqslant j \leqslant i\leqslant n.
\end{equation}
Therefore, we obtain a Gelfand-Tsetlin pattern of shape $\lambda$:
\begin{equation}
	\begin{array}{ccccccc}
		\tau_{1}&~&\tau_{2}&~&...&~&\tau_{n+1}\\
		~&\tau_{n,1}&~&...&~&\tau_{n,n}&~\\
		~&~&...&...&...&~&~\\
		~&~&\tau_{2,1}&~&\tau_{2,2}&~&~\\
		~&~&~&\tau_{1,1}&~&~&~
	\end{array}
\end{equation}
So we have a one-to-one correspondence between $\mathcal{Q}(\lmd)$ and the Gelfand-Tsetlin basis. Note that there is one-to-one correspondence between $\mathcal{B}(\lmd)$ and $\mathcal{Q}(\lmd)$. In summary, we have:\pse

{\bf Theorem 3.1}\quad {\it The set $\mathcal{B}(\lambda)$ is a basis of $V_n(\lambda)$.  Moreover,
	there is a bijection between $\mathcal{B}(\lambda)$ and the Gelfand-Tsetlin basis.}\psp

Next we want to find all the $sl(n)$ singular vectors in $\wht M$ and $\ol M=V_n(\lmd)$ by the method of determinants given in \cite{X4}.	Set
	\begin{equation}
		y_{i,j}=\left\lbrace
		\begin{aligned}
			&D_{n,i}(i), &\mathrm{if}~ i=j, \\
			&x_{i,j}, &\mathrm{if}~ 1\leqslant j<i \leqslant n
		\end{aligned}
		\right.
	\end{equation}
	Indeed,
	\begin{eqnarray}
			y_{k,k}&=&\begin{vmatrix}
				x_{k,k}&-1&0&...&0\\
				x_{k+1,k}&x_{k+1,k+1}&-1&...&0\\
				...&...&...&...&...\\
				x_{n-1,k}&x_{n-1,k+1}&x_{n-1,k+2}&...&-1\\
				x_{n,k}&x_{n,k+1}&x_{n,k+2}&...&x_{n,n}
			\end{vmatrix}\nonumber\\
			&=&x_{k,k}y_{k+1,k+1}+\begin{vmatrix}
				x_{k+1,k}&-1&0&...&0\\
				x_{k+2,k}&x_{k+2,k+2}&-1&...&0\\
				...&...&...&...&...\\
				x_{n-1,k}&x_{n-1,k+2}&x_{n-1,k+3}&...&-1\\
				x_{n,k}&x_{n,k+2}&x_{n.k+3}&...&x_{n,n}\\
			\end{vmatrix}.	
	\end{eqnarray}
		Thus
\begin{equation}x_{k,k}=y_{k+1,k+1}^{-1}y_{k,k}+\sum\limits_{l=k+1}^{n}g_{l}y_{l,k}, \end{equation}where
\begin{equation}g_{l}\mbox{'s are rational functions in}\; \left\lbrace y_{i,j} |k+1\leqslant j\leqslant i\leqslant n\right\rbrace. \end{equation}
So the sets of variables
\begin{equation}	 \mathcal{X}=\left\lbrace x_{i,j} | 1\leqslant j\leqslant i \leqslant n\right\rbrace \;\mbox{and}\;\mathcal{Y}=\left\lbrace y_{i,j} | 1\leqslant j\leqslant i \leqslant n\right\rbrace \end{equation}
	are functionally equivalent. Let $f$ be a rational function in $\mathcal{X}$, which is annihilated by positive root vectors $\{E_{s,t}\mid 1\leq t<s\leq n\}$. Then it can also be written as a rational function in $\mathcal{Y}$. We prove that $f$ is independent of $\left\lbrace y_{i,j} | 1\leqslant j< i \leqslant n\right\rbrace$ by backward induction. By (2.61) and (3.19), \begin{equation}E_{s,t}(y_{k,k})=0\qquad\for\;\mbox{any}\;\;1\leqslant k\leqslant n\;\;\mbox{and}\;\;1\leqslant t<s \leqslant n.\end{equation}
In particular,
	\begin{equation}
E_{n,n-1}(f)=\sum\limits_{1\leqslant j< i \leqslant n}(x_{n,n}\ptl_{n,n-1}-\ptl_{x_{n-1,n-1}})(y_{i,j})\frac{\partial f}{\partial y_{i,j}}
=y_{n,n}\frac{\partial f}{\partial y_{n,n-1}}=0
	\end{equation}
by (2.61) and (3.18).
	Thus $f$ is independent of $y_{n,n-1}$. Assume that $f$ is independent of $\left\lbrace y_{i,j} | k\leqslant j< i \leqslant n\right\rbrace $ for some $ 1<k\leqslant n-1$. By (2.61), (3.20) and (3.23),
	\begin{eqnarray}	
			E_{k,k-1}(f)&=&\sum\limits_{1\leqslant j<i\leqslant n}E_{k,k-1}(y_{i,j})\frac{\partial f}{\partial y_{i,j}}
			=\sum\limits_{l=k+1}^{n}y_{l,k}\frac{\partial f}{\partial y_{l,k-1}}+x_{k,k}\frac{\partial f}{\partial y_{k,k-1}}\nonumber\\ &=&\sum\limits_{l=k+1}^{n}y_{l,k}\big( \frac{\partial f}{\partial y_{l,k-1}}+g_{l}\frac{\partial f}{\partial y_{k,k-1}}\big) +y_{k+1,k+1}^{-1}y_{k,k}\frac{\partial f}{\partial y_{k,k-1}}=0.
	\end{eqnarray}
	By inductive assumption, $f$ is independent of $\left\lbrace y_{l,k} | k+1 \leqslant l\leqslant n \right\rbrace $.
 Thus the above equation yields\begin{equation}
		y_{k+1,k+1}^{-1}y_{k,k}\frac{\partial f}{\partial y_{k,k-1}}=0.
	\end{equation}
	So $f$ is independent of $y_{k,k-1}$. Then (3.25) implies
	\begin{equation}
		\sum\limits_{l=k+1}^ny_{l,k}\frac{\partial f}{\partial y_{l,k-1}}=0\lra
\frac{\partial f}{\partial y_{k+1,k-1}}=\frac{\partial f}{\partial y_{k+1,k-1}}=\cdots=\frac{\partial f}{\partial y_{n,k-1}}=0.,
	\end{equation}
	which implies that $f$ is independent of $\left\lbrace y_{i,j} | k-1\leqslant j< i \leqslant n\right\rbrace $.
By induction, $f$ is independent of $\left\lbrace y_{i,j} | 1\leqslant j< i \leqslant n\right\rbrace $. Namely,
$f$ is  a rational function in $\{y_{1,1},y_{2,2},...,y_{n,n}\}$. According to (3.18), we have:\pse

 {\bf Lemma 3.2} \quad {\it A rational function in $\mathcal{X}$ annihilated by positive root vectors $\{E_{s,t}\mid 1\leq t<s\leq n\}$ must be a rational function in $\{D_{n,1}(1),D_{n,2}(2),...,D_{n,n}(n)\}$.}
\psp

 Suppose a rational function in $\left\lbrace D_{n,i}(i) | 1\leqslant i\leqslant n\right\rbrace $:
\begin{equation}
	\frac{F(D_{n,1}(1),...,D_{n,n}(n))}{G(D_{n,1}(1),...,D_{n,n}(n))}=H\in\mathscr{C}_n,
\end{equation}
 where $F$ and $G$ are polynomials in $n$ variables. Write
 \begin{equation}F(D_{n,1}(1),...,D_{n,n}(n)), G(D_{n,1}(1),...,D_{n,n}(n))\;\mbox{and}\;H\end{equation}
as polynomials in $\left\lbrace x_{i,j}~|~1\leqslant j\leqslant i\leqslant n-1\right\rbrace$ with coefficients in $\mathbb{C}[x_{n,1},...,x_{n,n}]$. Denote  the ``constant terms" by $f$, $g$ and $h$ respectively. Then
   \begin{equation}
   f(x_{n,1},...,x_{n,n})=g(x_{n,1},...,x_{n,n})h(x_{n,1},...,x_{n,n}).\label{fgh}
   \end{equation}
  Note that the ``constant term" in $D_{n,k}(k)$ is just $x_{n,k}$. Thus we have
\begin{equation}
	f(x_{n,1},...,x_{n,n})=F(x_{n,1},...,x_{n,n}),\;\;
	g(x_{n,1},...,x_{n,n})=G(x_{n,1},...,x_{n,n}).
\end{equation}
So (\ref{fgh}) becomes
\begin{equation}F(x_{n,1},...,x_{n,n})=G(x_{n,1},...,x_{n,n})h(x_{n,1},...,x_{n,n}). \end{equation}
 Substituting $D_{n,k}(k)$ into $x_{n,k}$ in the above equation, we get
 \begin{equation}\frac{F(D_{n,1}(1),...,D_{n,n}(n)}{G(D_{n,1}(1),...,D_{n,n}(n))}=h(D_{n,1}(1),...,D_{n,n}(n))\end{equation}
 is a polynomial of $\left\lbrace D_{n,i}(i) | 1\leqslant i\leqslant n\right\rbrace $. Therefore, every polynomial in $\mathscr{C}_n$ annihilated by $\left\lbrace E_{s,t}~|~1\leqslant t<s\leqslant n\right\rbrace $ must be a polynomial of $\left\lbrace D_{n,i}(i) | 1\leqslant i\leqslant n \right\rbrace $.

For $\be=(\be_1,\be_2,...,\be_n)\in\mbb{N}^n$, we set
	\begin{equation}\lmd(\be)=\sum\limits_{i=1}^{n-1}(k_{i+1}+\beta_i-\beta_{i+1})\lambda_{i}\end{equation}and
\begin{equation}\xi(\be)=(D_{n,n}(n))^{\beta_1}(D_{n,n-1}(n-1))^{\beta_{2}}\cdots(D_{n,1}(1))^{\beta_n}.\end{equation}
Moreover, we calculate that
\begin{equation}\mbox{the weight of}\;\;\xi(\be)=
\lmd(\be)\end{equation}by (2.1), (2.8) and (2.60).

Let $|\beta|$ be the length of the multi-index $\beta$, namely $|\beta|=\sum\limits_{i=1}^n\beta_i$. Set
\begin{equation}
		\Omega_r=\big\lbrace \be=( \beta_1,\beta_2,...,\beta_n)\in\mathbb{N}^{n}~|~|\beta|=r,~0\leqslant \beta_i\leqslant k_i,~i\in\overline{1,n} \big\rbrace,
	\end{equation}
\begin{equation}\Omega=\bigcup_{r=0}^{|\lmd|}\Omega_r
	\end{equation}
and
\begin{equation}
		\wht\Omega_r=\big\lbrace \be=(\beta_1,\beta_2,...,\beta_n)\in\mathbb{N}^{n}~|~|\beta|=r,~0\leqslant \beta_i\leqslant k_i,~i\in\overline{2,n} \big\rbrace.
	\end{equation}Then $\Omega_r\subset \wht\Omega_r$. Recall that a singular vector is a weight vector annihilated by positive root vectors. In the following, we identify a singular vector with its nonzero constant multiple.  Denote by $\msr A_r$ the subspace of homogeneous polynomials with degree $r$ in $\mbb C[x_{n,1},x_{n,2},...,x_{n,n}]$.
 By (2.98), (2.99), (3.35) and Lemma 3.2, we have:\pse

{\bf Theorem 3.3}\quad {\it For $r\in\mbb N$, the $sl(n)$ singular vectors in $\wht M_r$ are
$\{\xi(\be)\mid\be\in\wht \Omega_r\}$,
which generate all the $sl(n)$ irreducible components in
\begin{equation}\wht M_r=\msr A_r\otimes_{\mbb C}V_{n-1}(\sum_{i=1}^{n-1}k_{i+1}\lmd_i).\end{equation}
Numerically, it implies the sum-product identity
\begin{equation}\sum_{\be\in \wht\Omega_r}d_{n-1}(\lmd(\be))=\binom{n+r-1}{n-1}d_{n-1}(\sum_{i=1}^{n-1}k_{i+1}\lmd_i)
\end{equation}
(cf. (1.17)).

For $r\in\ol{0,|\lmd|}$, the $sl(n)$ singular vectors in $(V_n(\lmd))_r=\ol M_r$ are
$\{\xi(\be)\mid\be\in\Omega_r\}$. Moreover, $\{\xi(\be)\mid\be\in\Omega\}$ are the highest-weight vectors of
all the $sl(n)$ irreducible components in $V_n(\lmd)=\ol M$. The numeric $sl(n+1)\downarrow sl(n)$ branching rule
 \begin{equation}\sum_{\be\in \Omega}d_{n-1}(\lmd(\be))=d_n(\lmd)\end{equation}
naturally follows (cf. (1.17)). If $r\leq k_1$, then $\wht \Omega_r=\Omega_r$, which yields
$(V_n(\lmd))_r=\wht M_r$.}\psp\pse

Now we want to use (3.41) to derive a formula for $(V_n(\lmd))_r$ and $V_n(\lmd)$. For $s\in\overline{1,n}$, let \begin{equation}r_s=r-\sum\limits_{j=1}^s(k_j+1).\end{equation} Set
\begin{equation}\widehat{\Omega}^s_{r_s}=\{ \beta\in\mathbb{N}^n|~|\beta|=r_s;~\beta_i\leqslant k_{i-1}~\mathrm{for}~i\in\overline{2,s};~\beta_i\leqslant k_i~\mathrm{for}~i\in\overline{s+1,n} \},\end{equation}
\begin{eqnarray}
	\widetilde{\Omega}^s_{r_s}&=&\{ \beta\in\mathbb{N}^n|~|\beta|=r_s;~\beta_i\leqslant k_{i-1}~\mathrm{for}~i\in\overline{2,s};\nonumber\\ & &\beta_{s+1}\leqslant k_s+k_{s+1}+1;~\beta_i\leqslant k_i~\mathrm{for}~i\in\overline{s+2,n} \}.	
\end{eqnarray}
Since $r\leqslant\sum\limits_{i=1}^nk_i$, we have $r_n<0$. So $\widehat{\Omega}^{n}_{r_{n}}=\widetilde{\Omega}^{n}_{r_{n}}=\emptyset$. Obviously $\widehat{\Omega}^s_{r_s}\subseteq\widetilde{\Omega}^s_{r_s}$ and
\begin{equation}
	\widetilde{\Omega}^s_{r_s}\setminus \widehat{\Omega}^s_{r_s}=(0,\cdots, \stackrel{s+1}{k_{s+1}+1},\cdots,0)+\widehat{\Omega}^{s+1}_{r_{s+1}}.
\end{equation}
By (3.40) and (3.41) with $\sum\limits_{i=1}^{n-1}k_{i+1}\lambda_i$ replaced by $\sum\limits_{i=1}^sk_i\lambda_i+\lambda_s+\sum\limits_{i=s}^{n-1}k_{i+1}\lambda_i$ and $r$ replaced by $r_s$, note that $\mathscr{A}_{r_s}=\emptyset$ and $\widetilde{\Omega}^s_{r_s}=\emptyset$ if $r_s<0$. We have
\begin{eqnarray}
	& &\dim (\mathscr{A}_{r_s}\otimes V_{n-1}(\sum\limits_{i=1}^sk_i\lambda_i+\lambda_s+\sum\limits_{i=s}^{n-1}k_{i+1}\lambda_i))\nonumber\\ &=&\sum\limits_{\beta\in\widetilde{\Omega}^s_{r_s}}d_{n-1}(\sum\limits_{i=1}^sk_i\lambda_i+
\lambda_s+\sum\limits_{i=s}^{n-1}k_{i+1}\lambda_i+\sum\limits_{i=1}^{n-1}(\beta_i-\beta_{i+1})\lambda_i).\label{0}
\end{eqnarray}
Hence, by (3.46) and (3.47),
\begin{eqnarray}
	& &\sum\limits_{\beta\in\widetilde{\Omega}^s_{r_s}\setminus \widehat{\Omega}^s_{r_s}}d_{n-1}(\sum\limits_{i=1}^sk_i\lambda_i+\lambda_s+\sum\limits_{i=s}^{n-1}k_{i+1}\lambda_i+\sum\limits_{i=1}^{n-1}(\beta_i-\beta_{i+1})\lambda_i)\nonumber\\
	&=&\sum\limits_{\beta\in\widehat{\Omega}^{s+1}_{r_{s+1}}}d_{n-1}(\sum\limits_{i=1}^sk_i\lambda_i+\lambda_s+\sum\limits_{i=s}^{n-1}k_{i+1}\lambda_i+\sum\limits_{i=1}^{n-1}(\beta_i-\beta_{i+1})\lambda_i+(k_{s+1}+1)(\lambda_{s+1}-\lambda_s))\nonumber\\
	&=&\sum\limits_{\beta\in\widehat{\Omega}^{s+1}_{r_{s+1}}}d_{n-1}(\sum\limits_{i=1}^{s+1}k_i\lambda_i+\lambda_{s+1}+\sum\limits_{i=s+1}^{n-1}k_{i+1}\lambda_i+\sum\limits_{i=1}^{n-1}(\beta_i-\beta_{i+1})\lambda_i)\nonumber\\
	&=&( \sum\limits_{\beta\in\widetilde{\Omega}^{s+1}_{r_{s+1}}}-\sum\limits_{\beta\in\widetilde{\Omega}^{s+1}_{r_{s+1}}\setminus \widehat{\Omega}^{s+1}_{r_{s+1}}})d_{n-1}(\sum\limits_{i=1}^{s+1}k_i\lambda_i+\lambda_{s+1}+\sum\limits_{i=s+1}^{n-1}k_{i+1}\lambda_i+\sum\limits_{i=1}^{n-1}(\beta_i-\beta_{i+1})\lambda_i)\nonumber\\
	&=&\dim (\mathscr{A}_{r_{s+1}}\otimes V_{n-1}(\sum\limits_{i=1}^{s+1}k_i\lambda_i+\lambda_{s+1}+\sum\limits_{i=s+1}^{n-1}k_{i+1}\lambda_i))\nonumber\\
	& &-\sum\limits_{\beta\in\widetilde{\Omega}^{s+1}_{r_{s+1}}\setminus \widehat{\Omega}^{s+1}_{r_{s+1}}}d_{n-1}(\sum\limits_{i=1}^{s+1}k_i\lambda_i+\lambda_{s+1}
+\sum\limits_{i=s+1}^{n-1}k_{i+1}\lambda_i+\sum\limits_{i=1}^{n-1}(\beta_i-\beta_{i+1})\lambda_i).
\end{eqnarray}
Moving items, we get
\begin{eqnarray}
	& &\dim (\mathscr{A}_{r_{s+1}}\otimes V_{n-1}(\sum\limits_{i=1}^{s+1}k_i\lambda_i+\lambda_{s+1}+\sum\limits_{i=s+1}^{n-1}k_{i+1}\lambda_i))\nonumber\\
	&=&\sum\limits_{\beta\in\widetilde{\Omega}^s_{r_s}\setminus \widehat{\Omega}^s_{r_s}}d_{n-1}(\sum\limits_{i=1}^sk_i\lambda_i+\lambda_s+\sum\limits_{i=s}^{n-1}k_{i+1}\lambda_i+\sum\limits_{i=1}^{n-1}(\beta_i-\beta_{i+1})\lambda_i)\nonumber\\
	& &+\sum\limits_{\beta\in\widetilde{\Omega}^{s+1}_{r_{s+1}}\setminus \widehat{\Omega}^{s+1}_{r_{s+1}}}d_{n-1}(\sum\limits_{i=1}^{s+1}k_i\lambda_i+\lambda_{s+1}
+\sum\limits_{i=s+1}^{n-1}k_{i+1}\lambda_i+\sum\limits_{i=1}^{n-1}(\beta_i-\beta_{i+1})\lambda_i).\label{1}
\end{eqnarray}
Note that $\widetilde{\Omega}^{n-1}_{r_{n-1}}\setminus \widehat{\Omega}^{n-1}_{r_{n-1}}=(0,...,0,(k_n+1))+\widehat{\Omega}^n_{r_n}=\varnothing$, multiply (\ref{1}) by $(-1)^s$ and sum over $1\leqslant s\leqslant n-2$,
\begin{eqnarray}
	& &\sum\limits_{s=1}^{n-2}(-1)^s\dim (\mathscr{A}_{r_{s+1}}\otimes V_{n-1}(\sum\limits_{i=1}^{s+1}k_i\lambda_i+\lambda_{s+1}+\sum\limits_{i=s+1}^{n-1}k_{i+1}\lambda_i))\nonumber\\
	&=&-\sum\limits_{\beta\in\widetilde{\Omega}^1_{r_1}\setminus \widehat{\Omega}^1_{r_1}}d_{n-1}((k_1+k_2+1)\lambda_1+\sum\limits_{i=2}^{n-1}k_{i+1}\lambda_i+\sum\limits_{i=1}^{n-1}(\beta_i-\beta_{i+1})\lambda_i)\nonumber\\
	& &+(-1)^{n-2}\sum\limits_{\beta\in\widetilde{\Omega}^{n-1}_{r_{n-1}}\setminus \widehat{\Omega}^{n-1}_{r_{n-1}}}d_{n-1}(\sum\limits_{i=1}^{n-2}k_i\lambda_i+(k_{n-1}+k_n+1)\lambda_{n-1}+\sum\limits_{i=1}^{n-1}(\beta_i-\beta_{i+1})\lambda_i)\nonumber\\
	&=&-\sum\limits_{\beta\in\widetilde{\Omega}^1_{r_1}\setminus \widehat{\Omega}^1_{r_1}}d_{n-1}((k_1+k_2+1)\lambda_1+\sum\limits_{i=2}^{n-1}k_{i+1}\lambda_i
+\sum\limits_{i=1}^{n-1}(\beta_i-\beta_{i+1})\lambda_i).\label{0.5}
\end{eqnarray}
Since $\widehat{\Omega}_r\setminus\Omega_r=(k_1+1,0,...,0)+\widehat{\Omega}^1_{r_1}$, as a conclusion, by (\ref{0}) and (\ref{0.5}),
\begin{eqnarray}
	& &\dim(\widehat{M}_r/(V_n(\lambda))_r)\nonumber\\
	&=&\sum\limits_{\beta\in\widehat{\Omega}^1_{r_1}}d_{n-1}((k_1+k_2+1+\beta_1-\beta_2)\lambda_1+\sum\limits_{i=1}^{n-1}(k_{i+1}+\beta_i-\beta_{i+1})\lambda_i)\nonumber\\
	&=&( \sum\limits_{\beta\in\widetilde{\Omega}^1_{r_1}}-\sum\limits_{\beta\in\widetilde{\Omega}^1_{r_1}\setminus\widehat{\Omega}^1_{r_1}}) d_{n-1}((k_1+k_2+1+\beta_1-\beta_2)\lambda_1+\sum\limits_{i=1}^{n-1}(k_{i+1}+\beta_i-\beta_{i+1})\lambda_i)\nonumber\\
	&=&\sum\limits_{s=1}^{n-1}(-1)^{s-1}\dim (\mathscr{A}_{r_s}\otimes V_{n-1}(\sum\limits_{i=1}^sk_i\lambda_i+\lambda_s+\sum\limits_{i=s}^{n-1}k_{i+1}\lambda_i))\nonumber\\
	&=&\sum\limits_{s=1}^{n-1}(-1)^{s-1}\binom{n+r-1-\sum\limits_{j=1}^s(k_j+1)}{n-1}d_{n-1}(\sum\limits_{i=1}^sk_i\lambda_i+\lambda_s+\sum\limits_{i=s}^{n-1}k_{i+1}\lambda_i)
\end{eqnarray}
Here we make the convention that for $p\in\mathbb{Z}$ and $q\in\mathbb{N}$,
\begin{equation}
	\binom{p}{q}=\left\lbrace
	\begin{aligned}
		&\frac{p!}{q!(p-q)!}, ~\mathrm{if}~p\geqslant q,\\
		&0, \qquad\qquad\mathrm{if}~p<q.
	\end{aligned}
\right.
\end{equation}

Moreover, we have
\begin{eqnarray}
	&&\dim (V_n(\lambda))_r=\dim\widehat{M}_r-\dim (\widehat{M}_r/(V_n(\lambda))_r)\nonumber\\
	&=&\binom{n+r-1}{n-1}d_{n-1}(\sum\limits_{i=1}^{n-1}k_{i+1}\lambda_i)\nonumber\\
	& &+\sum\limits_{s=1}^{n-1}(-1)^{s}\binom{n+r-1-\sum\limits_{j=1}^s(k_j+1)}{n-1}d_{n-1}(\sum\limits_{i=1}^sk_i\lambda_i+\lambda_s+\sum\limits_{i=s}^{n-1}k_{i+1}\lambda_i)\nonumber\\
	&=&\sum\limits_{s=0}^{n-1}(-1)^{s}\binom{n+r-1-\sum\limits_{j=1}^s(k_j+1)}{n-1}d_{n-1}(\sum\limits_{i=1}^sk_i\lambda_i+\lambda_s+\sum\limits_{i=s}^{n-1}k_{i+1}\lambda_i)
\end{eqnarray}
and
\begin{eqnarray}
	\dim V_n(\lambda)&=&\sum\limits_{r=1}^{|\lambda|}\dim (V_n(\lambda))_r\nonumber\\	&=&\sum\limits_{r=1}^{|\lambda|}\sum\limits_{s=0}^{n-1}(-1)^{s}\binom{n+r-1-\sum\limits_{j=1}^s(k_j+1)}{n-1}d_{n-1}(\sum\limits_{i=1}^sk_i\lambda_i+\lambda_s+\sum\limits_{i=s}^{n-1}k_{i+1}\lambda_i)\nonumber\\
&=&\sum\limits_{s=0}^{n-1}(-1)^{s}\binom{n+|\lambda|-\sum\limits_{j=1}^s(k_j+1)}{n}d_{n-1}(\sum\limits_{i=1}^sk_i\lambda_i+\lambda_s+\sum\limits_{i=s}^{n-1}k_{i+1}\lambda_i)\nonumber\\
	&=&\sum\limits_{s=0}^{n-1}(-1)^{s}\binom{\sum\limits_{j=s+1}^n(k_j+1)}{n}d_{n-1}
(\sum\limits_{i=1}^sk_i\lambda_i+\lambda_s+\sum\limits_{i=s}^{n-1}k_{i+1}\lambda_i) ,
\end{eqnarray}
where we just treat $\lambda_0=0$. \psp

{\bf Corollary 3.4}\quad {\it We have
\begin{equation} \dim (V_n(\lambda))_r=\sum\limits_{s=0}^{n-1}(-1)^{s}\binom{n+r-1-\sum\limits_{j=1}^s(k_j+1)}{n-1}d_{n-1}
(\sum\limits_{i=1}^sk_i\lambda_i+\lambda_s+\sum\limits_{i=s}^{n-1}k_{i+1}\lambda_i)\label{2}\end{equation}
for $r\in\ol{0,|\lmd|}$. Moreover,
\begin{equation}\sum\limits_{s=0}^{n-1}(-1)^{s}\binom{\sum\limits_{j=s+1}^n(k_j+1)}{n}d_{n-1}
(\sum\limits_{i=1}^sk_i\lambda_i+\lambda_s+\sum\limits_{i=s}^{n-1}k_{i+1}\lambda_i)=d_n(\sum_{\ell=1}^nk_\ell\lmd_\ell). \label{3}
\end{equation}}\psp\pse

Consider the case $k_1=k_2=\cdots=k_n=k$. The $sl(n+1)$-module $V_n(k\sum_{i=1}^n\lmd_i)$ is called a {\it Steinberg module}.
According to (1.17), its dimension
\begin{equation} d_n(k\sum_{i=1}^n\lmd_i)=(k+1)^{\frac{n(n+1)}{2}}.\end{equation}
Moreover,
\begin{equation}
	d_{n-1}(\sum\limits_{i=1}^sk\lambda_i+\lambda_s+\sum\limits_{i=s}^{n-1}k\lambda_i)=\binom{n}{s}(k+1)^{\frac{n(n-1)}{2}}
\end{equation}
By (\ref{2}), the dimension of homogeneous space of degree $r\in\overline{0,nk}$ is
\begin{eqnarray}
	\dim (V_n(\lambda))_r&=&\sum\limits_{s=0}^{n-1}(-1)^{s}\binom{n+r-1-s(k+1)}{n-1}d_{n-1}(\sum\limits_{i=1}^sk\lambda_i+\lambda_s+\sum\limits_{i=s}^{n-1}k\lambda_i)\nonumber\\
	&=&\sum\limits_{s=0}^{n-1}(-1)^{s}\binom{n+r-1-s(k+1)}{n-1}\binom{n}{s}(k+1)^{\frac{n(n-1)}{2}}
\end{eqnarray}
and (\ref{3}) becomes
\begin{eqnarray}	(k+1)^{\frac{n(n+1)}{2}}&=&\sum\limits_{s=0}^{n-1}(-1)^{s}\binom{(n-s)(k+1)}{n}d_{n-1}(\sum\limits_{i=1}^sk\lambda_i+\lambda_s+\sum\limits_{i=s}^{n-1}k\lambda_i)\nonumber\\
	&=&\sum\limits_{s=0}^{n-1}(-1)^{s}\binom{(n-s)(k+1)}{n}\binom{n}{s}(k+1)^{\frac{n(n-1)}{2}},
\end{eqnarray}
which is equivalent  to the well-known classical combinatorial identity:
\begin{equation}
	\sum\limits_{s=0}^n(-1)^{s}\binom{(n-s)(k+1)}{n}\binom{n}{s}=(k+1)^n
\end{equation}
(e.g., cf. Page 51 in \cite{Ri}).\pse
\psp

{\bf Remark 3.5}\quad  By (2.60) with $i=1$, $\{\lmd(\be)\mid\be\in\Omega\}$ are distinct weights of $gl(n)$.
Thus the multiplicities of the $gl(n)$ irreducible components in $V_n(\lmd)$ are one; that is, the multiplicity-one theorem of $gl(n+1)\downarrow gl(n)$ holds.

By (3.37), $|\Omega_r|=|\Omega_{|\lmd|-r}|$. So $sl(n)$-modules $\ol{M}_r$ and $\ol{M}_{|\lmd|-r}$ have the same number of irreducible components. Moreover, (3.34), (3.37) and (3.39) give
\begin{equation}\{\lmd(\be)\mid \be\in\widehat\Omega_{k_1+1}\setminus \Omega_{k_1+1}\}=\{(k_1+1)\lmd_1+\sum_{i=1}^{n-1}k_{i+1}\lmd_i\}.\end{equation}
Thus
\begin{equation}\wht{M}_{k_1+1}/\ol{M}_{k_1+1}=V_{n-1}((k_1+1)\lmd_1+\sum_{i=1}^{n-1}k_{i+1}\lmd_i).\end{equation}
According to (3.34) and (3.37),\
\begin{equation}\{\lmd(\be)\mid \be\in\Omega_{|\lmd|}\}=\{\sum_{i=1}^{n-1}k_i\lmd_i\}.\end{equation}
Hence
\begin{equation}\ol{M}_{|\lmd|}=V_{n-1}(\sum_{i=1}^{n-1}k_i\lmd_i)\end{equation}
is an irreducible $sl(n)$-module.

Note that as an $sl(n)$-module, $\msr A_r=V_{n-1}(r\lmd_1)$ in (3.40). By using the $sl(n)$ singular vectors $\{\xi(\be)\mid\be\in\wht \Omega_r\}$
in
\begin{equation}\wht M_r=V_{n-1}(r\lmd_1)\otimes V_{n-1}(\sum_{i=1}^{n-1}k_{i+1}\lmd_i),\end{equation}
one can derive the corresponding Clebsch-Gordan coefficients. Moreover, they can also be used to construct
 $sl(n)$-module homomorphisms
 \begin{equation}\mbox{from}\;V_{n-1}(r\lmd_1)\otimes_{\mbb C} V_{n-1}(\sum_{i=1}^{n-1}k_{i+1}\lmd_i)
 \;\mbox{to}\; V_{n-1}(\mu)\end{equation}for any dominant $sl(n)$ integral weight $\mu$,
 which are fundamental in constructing corresponding intertwining operators of the vertex operator algebra associated with the affine Kac-Moody Lie algebra $A^{(1)}$ (e.g., cf. (\cite{X1})). Using these intertwining operators, one can obtain the corresponding  exact solutions of Knizhnik-Zamolodchikov equation in WZW model of conformal field theory (cf. \cite{KZ, TK, X1}).

\section{Special Cases}

In this section, we give more detailed description on the special cases when the highest weights are
$k_1\lmd_1+k_2\lmd_2, k_1\lmd_1+k_2\lmd_n$ and $k\lmd_i$, respectively.

\subsection{$\lambda=k_1\lambda_{1}+k_2\lambda_{2}$, $n\geqslant 3$}

    In this case, (2.99) becomes
   \begin{eqnarray}
	S_n(\lambda)&=&\big\{
\prod\limits_{n-1\leqslant j\leqslant i\leqslant n}\prod\limits_{\Theta\in \mbb J_{i,j}} d_{i,j}(\Theta)^{\alpha_{i,j}(\Theta)} \mid \alpha_{i,j}(\Theta)\in\mbb{N};\nonumber\\ & &
\sum\limits_{i=j}^n\sum\limits_{\Theta\in\mbb J_{i,j}}\alpha_{i,j}(\Theta)\leqslant k_{n-j+1}~\mathrm{for}~j\in\overline{n-1,n}\big\}
\end{eqnarray}
because $k_3=k_4=\cdots=k_n=0$. According to (2.48), (4.1) only involves variables $\{x_{n-1,i},x_{n,j}\mid i\in\ol{1,n-1},
j\in\ol{1,n}\}$. Since $V_n(\lmd)$ is spanned by $S_n(\lmd)$, we have
\begin{equation}V_n(\lmd)\subset\mbb C[x_{n-1,i},x_{n,j}\mid i\in\ol{1,n-1},j\in\ol{1,n}].\end{equation} For simplicity, we redenote
\begin{equation}y_i=x_{n-1,i},\;\;x_j=x_{n,j}\qquad\for\;\; i\in\ol{1,n-1},
j\in\ol{1,n}.\end{equation}We also take the convention
\begin{equation}y_n=-1.\end{equation}Then the set
\begin{eqnarray}&&\big\{\prod\limits_{i=1}^nx_i^{\alpha_i}\prod\limits_{1\leqslant s<t\leqslant n}(x_sy_t-x_ty_s)^{\gamma_{s,t}}\prod\limits_{j=1}^{n-1}y_j^{\beta_j}\mid \al_i,\be_j,\gm_{s,t}\in\mbb N;\nonumber
\\ & &\sum_{i=1}^n\al_i\leq k_1;\sum_{1\leqslant s<t\leqslant n}\gm_{s,t}+\sum_{j=1}^{n-1}\be_j\leq k_2\big\}
\end{eqnarray}
spans $V_n(\lmd)$ by (4.1).

 According (2.60)-(2.62), we have the following representation formulas of $sl(n+1)$:
\begin{equation}E_{i,j}|_{V_n(\lambda)}=x_i\partial_{x_{j}}+y_i\partial_{y_{j}},~~E_{i,n}|_{V_n(\lambda)}=x_i\partial_{x_{n}}
+y_{i}\big( \sum\limits_{s=1}^{n-1}y_s\partial_{y_s}-k_2 \big),\end{equation}
\begin{equation}E_{i,n+1}|_{V_n(\lambda)}=x_i\big(\sum\limits_{s=1}^nx_s\partial_{x_s}-k_1 \big)+y_i\big( \sum\limits_{s=1}^{n-1}(x_s+x_ny_s)\partial_{y_s}\big)-k_2(x_i+x_ny_i),\end{equation}
\begin{equation}E_{n,j}|_{V_n(\lambda)}=x_n\partial_{x_{j}}-\partial_{y_j},~~E_{n+1,j}|_{V_n(\lambda)}=-\partial_{x_j},\end{equation}
\begin{equation}(E_{i,i}-E_{n+1,n+1})|_{V_n(\lambda)}=\sum\limits_{s=1}^nx_s\partial_{x_s}+x_i\partial_{x_{i}}+y_i\partial_{y_i}-k_1-k_2,\end{equation}
\begin{equation}(E_{n,n}-E_{n+1,n+1})|_{V_n(\lambda)}=\sum\limits_{s=1}^nx_s\partial_{x_s}+x_n\partial_{x_{n}}-\sum\limits_{s=1}^{n-1}y_s\partial_{y_s}-k_1.\end{equation}
\begin{equation}            	E_{n,n+1}|_{V_n(\lambda)}=x_n(\sum\limits_{s=1}^nx_s\partial_{x_s}-k_1)-\sum\limits_{s=1}^{n-1}(x_s+x_ny_s)\partial_{y_s},\;\;
E_{n+1,n}|_{V_n(\lambda)}
=-\partial_{x_n}
              \end{equation}
for $i,j\in\overline{1,n-1}$.

    By Theorem 3.3, $sl(n)$ singular vectors in $\wht M_r$ are
    \begin{equation}\{x_n^s(x_{n-1}+x_ny_{n-1})^t\mid s,t\in\mbb N; t\leq k_2; s+t=r\},\end{equation}whose weights are
     \begin{equation}\{(k_2+s-t)\lambda_1+t\lambda_2\mid s,t\in\mbb N; t\leq k_2; s+t=r\}\end{equation}for $r\in\mbb N$.
According to (1.17),
\begin{eqnarray}& &
    	d_{n-1}(\ell_1\lmd_1+\ell_2\lmd_2)\nonumber\\&=&\frac{\ell_1+1}{n-1}\binom{\ell_1+\ell_{2}+n-1}{n-2} \binom{\ell_{2}+n-2}{n-2}.\end{eqnarray}
    So (3.41) becomes
    \begin{eqnarray}& &\sum_{ s=\max\{0,r-k_2\}}^r\frac{k_2+2s-r+1}{n-1}\binom{k_{2}+s+n-1}{n-2} \binom{r-s+n-2}{n-2}
    \nonumber\\&=&\binom{n+r-1}{n-1}\binom{n+k_2-1}{n-1}.\end{eqnarray}

Let $k_3=\cdots=k_n=0$ in Corollary 3.4. Since $r\leqslant k_1+k_2<k_1+k_2+2$,
\begin{eqnarray}
	\binom{n+r-1-\sum\limits_{j=1}^s(k_j+1)}{n-1}=0,~\binom{\sum\limits_{j=s+1}^n(k_j+1)}{n}=\binom{n-s}{n}=0
\end{eqnarray}
when $s>1$. Thus (\ref{2}) and (\ref{3}) become
\begin{eqnarray} \dim(V_n(\lambda))_r&=&\binom{n+r-1}{n-1}d_{n-1}(k_2\lambda_1)-\binom{n+r-2-k_1}{n-1}d_{n-1}((k_1+k_2+1)\lambda_1)\nonumber\\
	&=&\binom{n+r-1}{n-1}\binom{k_2+n-1}{n-1}-\binom{n+r-2-k_1}{n-1}\binom{k_1+k_2+n}{n-1}
\end{eqnarray}and	
\begin{eqnarray}
	\dim V_n(\lambda)&=&\binom{k_1+k_2+n}{n}d_{n-1}(k_2\lambda_1)-\binom{k_2+n-1}{n}d_{n-1}((k_1+k_2+1)\lambda_1)\nonumber\\
	&=&\binom{k_1+k_2+n}{n}\binom{k_2+n-1}{n-1}-\binom{k_2+n-1}{n}\binom{k_1+k_2+n}{n-1}.
\end{eqnarray}

    \subsection{$\lambda=k_1\lambda_1+k_n\lambda_n$}

    In this case, the representation of $sl(n+1)$ is given in (2.60)-(2.62) with $k_2=\cdots=k_{n-1}=0$.
          Recall that
    \begin{equation}
    	D_{j,1}(1)=\begin{vmatrix}
    		x_{1,1}&-1&0&...&0\\
    		x_{2,1}&x_{2,2}&-1&...&0\\
    		...&...&...&...&...\\
    		x_{j-1,1}&x_{j-1,2}&x_{j-1,3}&...&-1\\
    		x_{j,1}&x_{j,2}&x_{j,3}&...&x_{j,j}\\
    	\end{vmatrix}.
    \end{equation}
   Now the $sl(n+1)$-module $V_n(\lambda)$ is spanned by
    \begin{equation}
    	S_n(\lambda)=\left\lbrace \prod\limits_{i=1}^nx_{n,i}^{\alpha_i}\prod\limits_{j=1}^n(D_{j,1}(1))^{\beta_j}~\Bigg| \al,\be\in\mbb N^n;\;|\al|\leq k_1,\;|\be|\leq k_n\right\rbrace.
    \end{equation}
     By Theorem 3.3, $sl(n)$ singular vectors in $\wht M_r$ are
    \begin{equation}\{x_n^s(D_{n,1}(1))^t\mid s,t\in\mbb N; t\leq k_n; s+t=r\},\end{equation}whose weights are
     \begin{equation}\{s\lambda_1+(k_n-t)\lmd_{n-1}\mid s,t\in\mbb N; t\leq k_n; s+t=r\}\end{equation}for $r\in\mbb N$. So (3.41) becomes
    \begin{eqnarray}& &\sum_{ s=\max\{0,r-k_n\}}^r\frac{k_n+2s-r+n-1}{n-1}\binom{s+n-2}{n-2}\binom{k_n+s-r+n-2}{n-2}
    \nonumber\\&=&\binom{n+r-1}{n-1}\binom{n+k_n-1}{n-1}.\end{eqnarray}

Suppose $k_1+k_n\geq r\geq k_1+1$. The set
\begin{equation}\big\lbrace \prod\limits_{i=1}^nx_{n,i}^{\alpha_i}\prod\limits_{j=1}^n(D_{j,1}(1))^{\beta_j}\mid\al,\be\in\mbb N^n;|\al|= k_1;\beta_n=r-k_1;|\be|-\be_n\leq k_1+k_n-r\big\rbrace\end{equation}
forms a basis of $(V_n(\lmd))_r$ by (4.19). Thus
\begin{equation}\dim (V_n(\lambda))_r=\binom{k_1+n-1}{n-1}\binom{k_1+k_n-r+n-1}{n-1}.\end{equation}
By Theorem 3.3, $sl(n)$ singular vectors in $(V_n(\lmd))_r$ are
    \begin{equation}\{x_n^s(D_{n,1}(1))^t\mid s,t\in\mbb N; s\leq k_1;t\leq k_n; s+t=r\},\end{equation}whose weights are
     \begin{equation}\{s\lambda_1+(k_n-t)\lmd_{n-1}\mid s,t\in\mbb N; s\leq k_1;t\leq k_n; s+t=r\}\end{equation}for $r\in\mbb N$. So we have
 \begin{eqnarray}&&\sum_{ s=\max\{0,r-k_n\}}^{\min\{k_1,r\}}\frac{k_n+2s-r+n-1}{n-1}\binom{s+n-2}{n-2}\binom{k_n+s-r+n-2}{n-2}
\nonumber\\ &=&\binom{k_1+n-1}{n-1}\binom{k_1+k_n-r+n-1}{n-1}.\end{eqnarray}

    \subsection{$\lambda=k\lambda_i,~2\leqslant i\leqslant n-1$}

   For simplicity, we denote
   \begin{equation}
    	d_j(\Theta)=d_j(\theta_1,\theta_2...,\theta_{i+j-n})=
    	\begin{vmatrix}
    		x_{n-i+1,\theta_1}&x_{n-i+1,\theta_2}&...&x_{n-i+1,\theta_{i+j-1}}\\
    		x_{n-i+2,\theta_1}&x_{n-i+2,\theta_2}&...&x_{n-i+2,\theta_{i+j-1}}\\
    		...&...&...&...\\
    		x_{j,\theta_1}&x_{j,\theta_2}&...&x_{j,\theta_{i+j-1}}\\
    	\end{vmatrix}
    \end{equation}
    for $j\in\overline{n-i+1,n}$ and $\Theta=(\theta_1,...,\theta_{i+j-1})\in\mathbb{J}_{i,j}$. Set
    \begin{equation}
    	\mathbb{I}_{i,j}=\left\lbrace(\theta_1,...,\theta_{i+j-1})\in\mathbb{N}^{i+j-1}~|~\theta_1\leqslant n-i+1,~\theta_1<\theta_2<....<\theta_{i+j-n}\leqslant j \right\rbrace.
    \end{equation}
By Theorem 3.3, the homogeneous subspace $(V_n(\lambda))_r$  is spanned by
    \begin{equation}
    \big\lbrace \prod\limits_{j=n-i+1}^n\prod\limits_{\Theta\in \mathbb{I}_{i,j}}d_j(\Theta)^{\alpha_j(\Theta)}\mid\alpha_j(\Theta)\in\mbb N;\sum\limits_{j=n-i+1}^n\sum\limits_{\Theta\in \mathbb{I}_{i,j}}\alpha_j(\Theta)\leqslant k;\sum\limits_{\Theta\in \mathbb{I}_{i,n}}\alpha_n(\Theta)=r\big\rbrace.
    \end{equation}
So $V_n(\lmd)$ dose not involve the variables $\{x_{r,s}\mid 1\leq s\leq r\leq n-i\}$. The  representation of $sl(n+1)$ is
 given in (2.60)-(2.62) with $k_r=0$ for $i\neq r\in\ol{1,n}$, $|_{\msr C_n}$ replaced by $|_{V_n(\lmd)}$ and all the ingredients containing $x_{r,s}$ deleted for $1\leq s\leq r\leq n-i$.

     By Theorem 3.3, $sl(n)$ singular vectors in $\wht M_r$ are
    \begin{equation}\{x_n^s(D_{n,n-i+1}(n-i+1))^t\mid s,t\in\mbb N; t\leq k; s+t=r\},\end{equation}whose weights are
     \begin{equation}\{s\lambda_1+(k-t)\lmd_{i-1}+t\lambda_i\mid s,t\in\mbb N; t\leq k; s+t=r\}\end{equation}for $r\in\mbb N$.
    Moreover,
     \begin{equation}
    	d_{n-1}(k\lambda_{i-1})=\prod\limits_{j=1}^{i-1}\frac{\binom{k+n-i+j}{k}}{\binom{k+j-1}{k}}.
    \end{equation}
         So (3.41) becomes
    \begin{eqnarray}& &\sum_{ s=\max\{0,r-k\}}^rd_{n-1}(s\lambda_1+(k+s-r)\lmd_{i-1}+(r-s)\lambda_i)
    \nonumber\\&=&\binom{n+r-1}{n-1}\prod\limits_{j=1}^{i-1}\frac{\binom{k+n-i+j}{k}}{\binom{k+j-1}{k}}.\end{eqnarray}

   Suppose $r\in\ol{0,k}$.  By Theorem 3.3, $sl(n)$-module $(V_n(\lmd))_r$ is an irreducible module with highest-weight vector $(D_{n,n-i+1}(n-i+1))^r$, whose weight is $(k-r)\lmd_{i-1}+r\lambda_i$. Thus
    \begin{equation}\dim(V_n(\lmd))_r=\binom{k-r+i-1}{i-1}\binom{r+n-i}{n-i}\prod\limits_{j=1}^{i-1}\frac{\binom{k+n-i+j}{k}}{\binom{k+j}{k}}. \end{equation}Hence (3.42) is directly equivalent to
    \begin{equation}
 \sum\limits_{r=0}^{k}\binom{k-r+i-1}{i-1}\binom{r+n-i}{n-i}=\binom{k+n}{k}.\end{equation}

 \end{document}